\documentclass{paper}

\usepackage{lineno,hyperref}
\modulolinenumbers[5]

\usepackage[tbtags]{amsmath}
\usepackage{amsfonts,amssymb,mathrsfs,amscd,comment}

\newtheorem{Lm}{Lemma}[section]
\newtheorem{corollary}{Corollary}[section]
\newtheorem{theorem}{Theorem}[section]

\newtheorem{Remark}{Remark}[section]
\newtheorem{Example}{Example}[section]
\newtheorem{Def}{Definition}[section]

\begin{document}


\title{On various approaches to Besov-type spaces of variable smoothness\thanks{This research was carried out with the financial support of the Russian Foundation for Basic Research (grant nos.
11-01-00744, 10-01-91331) and the program ``Development of the Scientific Potential of Higher Learning Institutions''
of the Ministry of Education and Science of the Russian Federation (project no.\ 2.1.1/1662).}}

\author{A.\,I.~Tyulenev}

\maketitle

\numberwithin{equation}{section}

\renewcommand{\theequation}{\arabic{section}.\arabic{equation}}

\begin{abstract}
The paper is concerned with Besov spaces  of variable smoothness $B^{\varphi_{0}}_{p,q}(\mathbb{R}^{n},\{t_{k}\})$, in which
the norms are defined in terms of convolutions with smooth functions. A~relation is found between the spaces $B^{\varphi_{0}}_{p,q}(\mathbb{R}^{n},\{t_{k}\})$
and the spaces $\widetilde{B}^{l}_{p,q,r}(\mathbb{R}^{n},\{t_{k}\})$, which were introduced earlier by the author.

\end{abstract}

%


\section{Introduction}
Besov-type spaces of variable smoothness have been studied extensively during the last 30 years.
We point out only a~few of the landmarks in the extensive literature on
Besov spaces  of variable smoothness (and their generalizations):  \cite{Ullrich}, \cite{Dachun}, \cite{Almeida}, \cite{Kempka}, \cite{Kempka2}, \cite{Be1}, \cite{Be2}.

For further purposes we shall introduce a special class of weight sequences. By a~weight sequence (which will be denoted by $\{s_{k}\}$, $\{t_{k}\}$)
we shall mean a~function sequence whose elements are measurable functions that are positive almost everywhere on~$\mathbb{R}^{n}$.

\begin{Def}
\label{Def1.1} {\rm
A weight sequence $\{s_{k}\}=\{s_{k}(\cdot)\}_{k=0}^{\infty}$ will be said to lie in $ Y^{\alpha_{3}}_{\alpha_{1},\alpha_{2}}$ if,
for $\alpha_{3}\geq 0$, $\alpha_{1},\alpha_{2} \in \mathbb{R}$,

\smallskip
1) $\frac{1}{C_{1}} 2^{\alpha_{1}(k-l)} \le \frac{s_{k}(x)}{s_{l}(x)} \le C_{1} 2^{\alpha_{2}(k-l)}$, $l \le k \in \mathbb{N}_{0}$, $x \in \mathbb{R}^{n}$;
\vskip-20pt
\begin{gather}
\noalign{}\label{1.1}
\end{gather}
\vskip-20pt
2) $ s_{k}(x) \le C_{2} s_{k}(y)(1+2^{k}|x-y|)^{\alpha_{3}}$, $k \in \mathbb{N}_{0}$, $x,y \in \mathbb{R}^{n}$;
}

\noindent here $C_1$, $C_2$ are positive and independent of $x, y$ and $k, l$.
\end{Def}

In the majority of papers known at present the Besov spaces of variable smoothness were defined in terms of the Fourier transform.

Let $B^{n}$ be the unit ball of $\mathbb{R}^{n}$, $\Psi_{0} \in S(\mathbb{R}^{n})$, $\Psi_{0}(x)=1$
for $x \in B^{n}$, $\operatorname{supp}\Psi_{0} \subset 2B^{n}$. For $j \in \mathbb{N}$, we set $\Psi_{j}(x):=\Psi_{0}(2^{-j}x)-\Psi_{0}(2^{-j+1}x)$, $x \in \mathbb{R}^{n}$.

\begin{Def}
\label{Def1.2}
Let $p,q \in (0,\infty]$, $\alpha_{1},\alpha_{2} \in \mathbb{R}$, $\alpha_{3} \geq 0$, $\{s_{k}\} \in Y^{\alpha_{3}}_{\alpha_{1},\alpha_{2}}$.
By $B^{\{s_{k}\}}_{p,q}(\mathbb{R}^{n})$ we shall denote the space of all distributions $f \in S'(\mathbb{R}^{n})$ with finite quasi-norm
\begin{equation}
\label{1.2}
\|f|B^{\{s_{k}\}}_{p,q}(\mathbb{R}^{n})\|:=\|s_{j}F^{-1}(\Psi_{j}F[f])|l_{q}(L_{p}(\mathbb{R}^{n}))\|.
\end{equation}
\end{Def}

In \eqref{1.2} the symbols $F$ and $F^{-1}$ denote, respectively, the direct and inverse Fourier transforms.

Definition \ref{Def1.2} is not satisfactory by the following reasons. First, a weight sequence $\{s_{k}\}$ contains the functions
which increase slowly at infinity. Moreover, each weight function $s_{k}$ has no points of singularity or degeneracy points.

An axiomatic approach to function spaces was developed in \cite{Ullrich}, \cite{Dachun}.
A more general (than~$L_{p}$) function space satisfying a~certain set of axioms was used on the right of~\eqref{1.2}.
This approach is capable of dealing with weighted Besov spaces of variable smoothness $\{s_{k}\} \in Y^{\alpha_{3}}_{\alpha_{1},\alpha_{2}}$
(with weights satisfying the doubling property).

Besov spaces of variable smoothness with norm defined in terms of classical differences were studied in \cite{Be1},  \cite{Be2} for $p,q \in (1,\infty)$.
Moreover, a~wider (in comparison with $Y^{\alpha_{3}}_{\alpha_{1},\alpha_{2}}$) class of weighted sequences $^{\text{\rm loc}}Y^{\alpha_{3}}_{\alpha_{1},\alpha_{2}}$
was examined---this class differs from the class $Y^{\alpha_{3}}_{\alpha_{1},\alpha_{2}}$ in that condition~2) in Definition \ref{Def1.1} is replaced by the condition

\smallskip
\noindent 2$'$)\ \ $s_{k}(x) \le C_{2} s_{k}(y)2^{k\alpha_{3}}$, $k \in \mathbb{N}_{0}$, $|x-y| \le 2^{-k}$.\hfill (1.3)
\smallskip

\setcounter{equation}{3}

We pursue two goals in this paper. First, we want to extend the theory of the spaces $B^{\{s_{k}\}}_{p,q}(\mathbb{R}^{n})$ with $p,q \in (0,\infty]$
to a~much wider class of weight sequences.
In doing so, the approach developed by Rychkov \cite{Ry} will be used instead of the Fourier analysis.
The class of weighted sequences to be introduced below turns out to be so broad that it will envelope both the weighted Besov spaces of~\cite{Ry} and the Besov spaces of variable smoothness
(only for constant $p$ and~$q$)
which were actively studied by H.~Kempka with coauthors in \cite{Kempka}, \cite{Kempka2} (for details see Remark~\ref{R2.4} below).
It is also worth noting that Besov spaces of variable smoothness are first dealt with in this degree of generality.

Second, we will establish the relation between the spaces appearing as a~result of such extension and
the spaces  $\widetilde{B}^{l}_{p,q,r}(\mathbb{R}^{n},\{t_{k}\})$, which were introduced by the author in~\cite{Tyu} using the methods of nonlinear spline approximation.
The space $\widetilde{B}^{l}_{p,p,r}(\mathbb{R}^{n},\{\gamma_{k}\})$ (with $1\le r <p < \infty$) arises in the natural way as the trace space of  the weighted Sobolev space
$W^{l}_{p}(\mathbb{R}^{n}, \gamma)$ with weight $\gamma \in A^{\rm loc}_{\frac{p}{r}}(\mathbb{R}^{n})$ (for details, see \S\,3 of~\cite{Tyu}).
Hence, the problem of finding an equivalent description of the space $\widetilde{B}^{l}_{p,q,r}(\mathbb{R}^{n},\{t_{k}\})$
is of great importance in some practical applications (for details see~\S\,5).

\section{Auxiliary results}

The following convention will be followed throughout the paper: when no limits of integration are explicitly shown,
it is to be understood that the integration is taken over the space~$\mathbb{R}^{n}$.

By $D(\mathbb{R}^{n})$ we shall denote the linear space of compactly supported infinitely differentiable functions;
the topology on this space is standard (see, for example, Ch.~6 of~\cite{Rudin}). Next, $D'(\mathbb{R}^{n})$ will denote the linear space of continuous linear functionals on
$D(\mathbb{R}^{n})$. Following~\cite{Ry} we shall be concerned only with functionals from the space
$S_{e}(\mathbb{R}^{n}) \subset D'(\mathbb{R}^{n})$. Note that $f \in S_{e}(\mathbb{R}^{n})$ if and only if
$$
|\langle f,\varphi\rangle| \le C_{f}\sup\{|D^{\alpha}\varphi(x)|\operatorname{exp}(N_{f}|x|),x \in \mathbb{R}^{n}, |\alpha| \le N_{f}\}, \quad \varphi  \in D(\mathbb{R}^{n}).
$$

Given a function $g \in D(\mathbb{R}^{n})$, we set $g_{j}:=2^{jn}g(\cdot 2^{j})$, $j \in \mathbb{N}$.

Next, for a function $\psi \in D(\mathbb{R}^{n})$ we let $L_{\psi}$ denote the supremum of~$L$ such that
\begin{equation}
\label{2.1}
\int x^{\beta}\psi(x)\,dx=0, \qquad  |\beta| \le L.
\end{equation}

A~function $\psi$ is said to have zero moment of order~$L$ if condition \eqref{2.1} is satisfied.

In what follows $Q^{n}$ will denote a~cube in the space $\mathbb R^n$  with sides parallel to coordinate axes,  $r(Q^{n})$
will denote the side length of~$Q^n$. For $\delta >0$, by $\delta Q^{n}$ we shall mean the cube, concentric with a~cube $Q^{n}$, with side length $r(\delta Q^{n}):=\delta r(Q^{n})$. By $Q^{n}_{k,m}$ we denote a~dyadic cube of rank~$k$. More precisely, $Q^{n}_{k,m}:=\prod\limits_{i=1}^{n}\left(\frac{m_{i}}{2^{k}},\frac{m_{i}+1}{2^{k}}\right)$, $\widetilde{Q}^{n}_{k,m}:=\prod\limits_{i=1}^{n}\left[\right.\frac{m_{i}}{2^{k}},\frac{m_{i}+1}{2^{k}}\left.\right)$ for $k \in \mathbb{N}_{0}$, $m \in \mathbb{Z}^{n}$. We set also $I^{n}:=(-1,1)^{n}$.

For $p \in (0,\infty]$ and a measurable set $E$ of positive measure, we  let $L_{p}(E)$ denote the set of all classes of equivalent measurable functions with finite quasi-norm $\|f|L_{p}(E)\|:=\left(\int\limits_{E}|f(x)|^{p}\,dx\right)^{\frac{1}{p}}$
(for $p=\infty$ the corresponding modifications are straightforward: the integral is replaced by $\operatorname{ess\,sup}$).
Given $q \in (0,\infty]$, we let $l_{q}$ denote the linear space of all real sequences $\{a_{k}\}$ with finite quasi-norm
$\|a_{k}|l_{q}\|:=\left(\sum\limits_{k=1}^{\infty}|a_{k}|^{q}\right)^{\frac{1}{q}}$ (for $q=\infty$ the corresponding norm is defined using the supremum).

By a weight we shall mean an arbitrary measurable function that is positive almost everywhere.
For the definition and basic properties of the weighted class $A^{\text{\rm loc}}_{p}(\mathbb{R}^{n})$ for $p \in (1,\infty]$ we refer the reader to~\cite{Ry}.

Throughout, by $c$ or $C$ we shall denote insignificant constants, which may be different in different inequalities.
The constants will not be labeled. Sometimes, when it is essential for the understanding of the exposition which follows,
we shall indicate the parameters on which some or other constant is dependent.

In the present paper we shall slightly modify the definition of the weight class $X^{\alpha_{3}}_{\alpha,\sigma,p}$, which was introduced by the author in~\cite{Tyu}.

A~sequence of weights $\{t_{k}\}$ is said to be $p$-admissible for $p \in (0,\infty]$ if  $t_{k} \in L^{\text{\rm loc}}_{p}(\mathbb{R}^{n})$, $k \in \mathbb{N}_{0}$.

Let $\{t_{k}\}$ be a $p$-admissible weight sequence. By $\{t_{k,m}\}$ we shall denote the multiple sequence defined by
\begin{equation}
\label{2.2}
t_{k,m} := \|t_{k}|L_{p}(Q^{n}_{k,m})\| , \qquad k \in \mathbb{N}_{0}, \ m \in \mathbb{Z}^{n},
\end{equation}
and by $\{\overline{t}_{k}\}$ we shall denote the weight sequence
$$
\overline{t}_{k}(x):=2^{\frac{kn}{p}}\sum\limits_{m \in \mathbb{Z}^{n}}t_{k,m}\chi_{\widetilde{Q}^{n}_{k,m}}(x), \qquad k \in \mathbb{N}_{0}, \ x \in \mathbb{R}^{n}.
$$

\begin{Def}
\label{Def2.1}
{\rm Let $p,\sigma_{1},\sigma_{2} \in (0,\infty]$, let $\alpha^{1}:=\{\alpha^{1}_{k}\}$, $\alpha^{2}:=\{\alpha^{2}_{k}\}$ be sequences of positive real numbers,
and let $\alpha_{3} \geq 0$. We set $\sigma:=(\sigma_{1},\sigma_{2})$, $\alpha:=( \alpha^{1},\alpha^{2})$. By $X^{\alpha_{3}}_{\alpha,\sigma,p}$ we shall denote the set of
 all  $p$-admissible weight sequences $\{t_{k}\}:=\{t_{k}\}_{k=0}^{\infty}$ such that, for some $C_{1},C_{2}>0$,

\smallskip
{\rm 1)}\ \ $\displaystyle \biggl(2^{kn}\int\limits_{Q^{n}_{k,m}}t^{p}_{k}(x)\biggr)^{\frac{1}{p}}\biggl(2^{kn}\int\limits_{Q^{n}_{k,m}}t^{-\sigma_{1}}_{j}(x)\biggr)^{\frac{1}{\sigma_{1}}} \le C_{1}\frac{\alpha^{1}_{k}}{\alpha^{1}_{j}}, \quad
 0 \le k \le j, \ m \in \mathbb{Z}^{n}$,\break\null\hfill (2.3)

{\rm 2)}\ \ $\displaystyle
\biggl(2^{kn}\int\limits_{Q^{n}_{k,m}}t^{p}_{k}(x)\biggr)^{-\frac{1}{p}}\biggl(2^{kn}\int\limits_{Q^{n}_{k,m}}t^{\sigma_{2}}_{j}(x)\biggr)^{\frac{1}{\sigma_{2}}} \le C_{2}\frac{\alpha^{2}_{j}}{\alpha^{2}_{k}}, \quad
 0 \le k \le j,\  m \in \mathbb{Z}^{n}$,\break\null\hfill (2.4)

\setcounter{equation}{4}
\smallskip

\noindent (for $p=\infty$ or $\sigma_{1}=\infty$ or $\sigma_{2}=\infty$ the modifications in (2.3) and (2.4) are standard).

{\rm 3)} for all  $k \in \mathbb{N}_{0}$
 \begin{equation}
\label{2.5}
0<t_{k,m} \le 2^{\alpha_{3}} t_{k,\widetilde{m}}, \qquad m,\widetilde{m} \in \mathbb{Z}^{n}, \ \ |m_{i}-\widetilde{m}_{i}| \le 1, \ \ i=1,\dots,n.
\end{equation}
}
\end{Def}

\begin{Remark}
\label{R2.1}
{\rm Let parameters $c_{1},c_{2} \geq 1$ and sequence $\{t_{k}\} \in X^{\alpha_{3}}_{\alpha,\sigma,p}$ be fixed.
Then it easily follows from \eqref{2.5} that

{\rm 1$'$)}\ $\displaystyle \biggl(2^{kn}\!\!\!\int\limits_{c_{1}Q^{n}_{k,m}}t^{p}_{k}(x)\biggr)^{\frac{1}{p}}\biggl(2^{kn}\int\limits_{c_{2}Q^{n}_{k,m}}t^{-\sigma_{1}}_{j}(x)\biggr)^{\frac{1}{\sigma_{1}}} \le C'_{1}\frac{\alpha^{1}_{k}}{\alpha^{1}_{j}}$,
$ 0 \le k \le j$, $ m \in \mathbb{Z}^{n}$,\break\null\hfill(2.6)

{\rm 2$'$)}\ \ $\displaystyle \biggl(2^{kn}\!\!\!
\int\limits_{c_{1}Q^{n}_{k,m}}t^{p}_{k}(x)\biggr)^{-\frac{1}{p}}\biggl(2^{kn}\int\limits_{c_{2}Q^{n}_{k,m}}t^{\sigma_{2}}_{j}(x)\biggr)^{\frac{1}{\sigma_{2}}} \le C'_{2}\frac{\alpha^{2}_{j}}{\alpha^{2}_{k}}$, $0 \le k \le j$, $m \in \mathbb{Z}^{n}$.\break\null\hfill(2.7)

\setcounter{equation}{7}

\noindent Here, the constants $C'_{1}, C'_{2}$ depend only on the constants $C_{1},C_{2},c_{1},c_{2},\alpha_{3},\sigma_{1},\sigma_{2}, p$.}
\end{Remark}

\begin{Remark}
\label{R2.2}
{\rm Let $\alpha_{1},\alpha_{2} \in \mathbb{R}$, $\{\alpha^{1}_{k}\}=\{2^{k\alpha_{1}}\}$, $\{\alpha^{2}_{k}\}=\{2^{k\alpha_{2}}\}$.
If under these assumptions one requires that conditions (2.3), (2.4) were satisfied for the weight sequence $\{\overline{t}_{k}\}$ (instead of $\{t_{k}\}$), then
Definition \ref{Def2.1} is equivalent to the corresponding definition from~\cite{Tyu}.}
\end{Remark}

The following elementary assertion will be useful in \S\,3.

\begin{Lm}
\label{Lm2.1}
Let $\alpha_{3} \geq 0$, $\sigma_{1},p \in (0,\infty]$, $\sigma_{2}=p$, let $\{\alpha^{1}_{k}\}$, $\{\alpha^{2}_{k}\}$
be sequences of positive real numbers. Next, let a~$p$-admissible weight sequence $\{t_{k}\} \in X^{\alpha_{3}}_{\alpha,\sigma,p}$. Then the weight sequence
$\{\overline{t}_{k}\} \in X^{\alpha_{3}}_{\alpha,\sigma,p}$.
\end{Lm}

\textbf{Proof.} It is clear that estimates (2.4), \eqref{2.5} holds with $\{\overline{t}_{k}\}$ instead of $\{t_{k}\}$.

Note that the exponent $\sigma_{1}$ may be written as  $\sigma_{1}=\theta p'_{\theta}$ with some  $\theta \in (0,p]$ (here we assume that
$\theta p'_{\theta}=\infty$ for $p'_{\theta}=\infty$ and $p_{\theta}=1$ for $p=\theta \in (0,\infty]$). Hence, by H\"older's inequality,
    \begin{equation}
\label{holder}
2^{-\frac{kn}{\theta}}=\biggl(\int\limits_{Q^{n}_{k,m}}\frac{t^{\theta}_{k}(x)}{t^{\theta}_{k}(x)}\,dx\biggr)^{\frac{1}{\theta}} \le \biggl(\int\limits_{Q^{n}_{k,m}}t^{p}_{k}(x)\,dx\biggr)^{\frac{1}{p}}\biggl(\int\limits_{Q^{n}_{k,m}}t^{-\sigma_{1}}_{k}(x)\,dx\biggr)^{\frac{1}{\sigma_{1}}}.
\end{equation}

A direct calculation using estimate \eqref{holder} gives (2.3) with $\{\overline{t}_{k}\}$ instead of $\{t_{k}\}$.

\begin{Remark}
\label{R2.3}
{\rm Let $p \in (0,\infty]$, $\alpha_{1},\alpha_{2} \in \mathbb{R}^{n}$, $\{\alpha^{1}_{k}\}=\{2^{k\alpha_{1}}\}$, $\{\alpha^{2}_{k}\}=\{2^{k\alpha_{2}}\}$.
The embedding $^{\text{\rm loc}}Y^{\alpha_{3}}_{\alpha_{1},\alpha_{2}} \subset X^{\alpha_{3}}_{\alpha,\sigma,p}$ for all  $\sigma_{1},\sigma_{2} \in (0,\infty]$
can be proved by elementary estimates. However, if, when considering  the class $X^{\alpha_{3}}_{\alpha,\sigma,p}$, one is content only with
such $p$-admissible weight sequences $\{t_{ k}\}$ for which $\{t_{k}\}=\{\overline{t}_{k}\}$, then
 $X^{\alpha_{3}}_{\alpha,\sigma,p} \subset ^{\text{\rm loc}}Y^{\alpha_{3}}_{\alpha_{1}-\frac{n}{\sigma_{1}},\alpha_{2}+\frac{n}{\sigma_{2}}}$.}
\end{Remark}

\begin{Example}
\label{Ex2.1}
{\rm We give an example illustrating the advantage of Definition \ref{Def2.1} over Definition \ref{Def1.1}.
Let $p \in (0,\infty)$, $r \in (0,p)$, a~weight  $\gamma^{p} \in A^{\text{\rm loc}}_{\frac{p}{r}}(\mathbb{R}^{n})$, and a~weight sequence
$\{s_{k}\} \in ^{\text{\rm loc}}Y^{\alpha_{3}}_{\alpha_{1},\alpha_{2}}$.
We set $t_{k}(x)=s_{k}(x)\gamma(x)$ with $x \in \mathbb{R}^{n}$, $k \in \mathbb{N}_{0}$.
Now a~direct calculation using the definition of the class $A^{\text{\rm loc}}_{\frac{p}{r}}(\mathbb{R}^{n})$ (see~\cite{Ry})
shows that $\{t_{k}\}\in X^{\alpha_{3}}_{\alpha,\sigma,p}$ with $\sigma_{1}=\frac{pr}{p-r}$, $\sigma_{2}=p$, $\{\alpha^{1}_{k}\}=\{2^{k\alpha_{1}}\}$,
$\{\alpha^{2}_{k}\}=\{2^{k\alpha_{2}}\}$ (with some $\alpha_{1},\alpha_{2} \in \mathbb{R}$).
 Thus, a~multiplication of a~fairly `good' sequence $\{s_{k}\}$ by a~sufficiently `bad' weight~$\gamma$ does not impair the exponents
 $\alpha_{1}$, $\alpha_{2}$ (of course, if we are dealing with the class $X^{\alpha_{3}}_{\alpha,\sigma,p}$).}
\end{Example}

\begin{Def}
\label{Def2.2}
{\rm Let $p,q,\sigma_{1},\sigma_{2} \in (0,\infty]$, $\alpha_{3} \geq 0$, and let $\alpha^{1}:=\{\alpha^{1}_{k}\}$, $\alpha^{2}:=\{\alpha^{2}_{k}\}$ be
sequences of positive numbers. Assume that a~weight sequence $\{t_{k}\} \in X^{\alpha_{3}}_{\alpha,\sigma,p}$. Let $\varphi_{0} \in C^{\infty}_{0}(\mathbb{R}^{n})$.
We set $\varphi(x):=\varphi_{0}(x)-2^{-n}\varphi_{0}(\frac{x}{2})$ for $x \in \mathbb{R}^{n}$. By $B^{\varphi_{0}}_{p,q}(\mathbb{R}^{n}, \{t_{k}\})$
we denote the set of all distributions  $f \in S'_{e}$ with finite quasi-norm
\begin{equation}
\label{Besov1}
\|f|B^{\varphi_{0}}_{p,q}(\mathbb{R}^{n}, \{t_{k}\})\|:= \biggl(\sum\limits_{k=0}^{\infty} \|t_{k} (\varphi_{k} \ast f)|L_{p}(\mathbb{R}^{n})\|^{q}\biggr)^{\frac{1}{q}}
\end{equation}
(the modifications of \eqref{Besov1} in the case $q=\infty$ are clear).
}
\end{Def}
Let $r \in (0,\infty]$, $f \in L^{loc}_{r}(\mathbb{R}^{n})$,  $l\in \mathbb N$. We set (the modification in the case $r=\infty$ is clear)
\begin{gather*}
  \Delta^l (h) f(x)= \sum^l_{i=0} (-1)^{l+i}C^i_l f(x+ih), \qquad x\in \mathbb R^n, \ h\in \mathbb R^n; \\
  \delta^l_r(x+2^{-k}I^n)f := \Bigl( 2^{2kn} \int\limits_{x+\frac{I^n}{2^k}} \int\limits_{\frac{I^n}{2^k}}
 |\Delta^l(h)f(y)|^r dh\,dy\Bigr)^{\frac{1}{r}}, \qquad x\in \mathbb R^n, \ k\in \mathbb N_0.
\end{gather*}

\begin{Def} [see \cite{Tyu}]
\label{Def2.3}
{\rm Let $p,q,r,\sigma_{1},\sigma_{2} \in (0,\infty]$, $\alpha_{3} \geq 0$, $\alpha^{1}:=\{\alpha^{1}_{k}\}$, $\alpha^{2}:=\{\alpha^{2}_{k}\}$ be sequences of
positive numbers. Assume that a~weight sequence $\{t_{k}\} \in X^{\alpha_{3}}_{\alpha,\sigma,p}$. By  $\widetilde{B}^{l}_{p,q,r}(\mathbb{R}^{n},\{t_{k}\})$ we denote the class of
functions $f \in L^{\text{\rm loc}}_{r}(\mathbb{R}^{n})$ with finite quasi-norm
\begin{multline}
\label{Besov2}
\widetilde{B}^{l}_{p,q,r}(\mathbb{R}^{n},\{t_{k}\}):= \\
:=\biggl(\sum\limits_{k=0}^{\infty} \|t_{k} \delta^{l}_{r}(\cdot+2^{-k}I^{n}) f|L_{p}(\mathbb{R}^{n})\|^{q}\biggr)^{\frac{1}{q}}+\|t_{0}\|f|L_{r}(\cdot+I^{n})\||L_{p}(\mathbb{R}^{n})\|
\end{multline}
(the modifications of \eqref{Besov2} in the case $q=\infty$ are clear).}
\end{Def}

\begin{Remark}
\label{R2.4}
{\rm Let $p \in (0,\infty)$, $s \in \mathbb{R}$, $\{\alpha^{1}_{k}\}=\{\alpha^{2}_{k}\}=\{2^{ks}\}$, a weight $\gamma^{p} \in A^{\text{\rm loc}}_{\infty}(\mathbb{R}^{n})$
and $L_\varphi \ge [s]$
Taking  $\{t_{k}\}=\{2^{ks}\gamma\}$ in the Definition~\ref{2.2} we obtain the definition of the weighted Besov space of~\cite{Ry}; we denote this space by  $B^{s}_{p,q}(\mathbb{R}^{n},\gamma)$.
It is known that in the case $\gamma \in A_{\infty}(\mathbb{R}^{n})$ and $L_{\varphi} \ge [s]$ the definition of the weighted Besov spaces, as given in terms of the
Fourier transform (see the introduction), and the definition in terms of convolutions with smooth functions are equivalent  (see~\cite{Ry} for details).

Arguing as in the proof of Proposition 1.2 of~\cite{Ry2} and using Theorems 3.6 and 3.8 of~\cite{Kempka3}, one may show that,
for $p,q \in (0,\infty]$, $\{s_{k}\} \in Y^{\alpha_{3}}_{\alpha_{1},\alpha_{2}}$ and $L_{\varphi} \geq [\alpha_{2}]$,  the space
$B^{\varphi_{0}}_{p,q}(\mathbb{R}^{n},\{s_{k}\})\cap S'(\mathbb R^n)$ coincides with the space $B^{\{s_{k}\}}_{p,q}(\mathbb{R}^{n})$ (Definition~\ref{Def1.2}),
the corresponding norms being equivalent.
}
\end{Remark}

\begin{theorem}[Hardy's inequality]
\label{Th2.1}
 Let $1 \le  s \le \infty$ and let $\{a_{k}\}$, $\{\beta_{k}\}$ be sequences of positive real numbers. Then the inequality
\begin{equation}
\label{2.10}
\biggl(\sum\limits_{k=0}^{\infty}\beta^{s}_{k}b_{k}^{s}\biggr)^{\frac{1}{s}} \le C \biggl(\sum\limits_{k=0}^{\infty}\beta^{s}_{k}a_{k}^{s}\biggr)^{\frac{1}{s}}
\end{equation}
{\rm  (}with straightforward modifications in the case $s=\infty)$
holds with constant $C > 0$ independent of the sequences $\{a_{k}\}$ in the cases:

\smallskip

{\rm 1)}\hfill
$\displaystyle
\begin{gathered}
b_{k} \le C' \biggl(\sum\limits_{j=k}^{\infty}a_{k}\biggr) \quad \text{if}\quad
\sup\limits_{n \in \mathbb{N}}\biggl(\sum\limits_{k=0}^{n}\beta^{s}_{k}\biggr)^{\frac{1}{s}}\biggl(\sum\limits_{k=n}^{\infty}\beta_{k}^{-s'}\biggr)^{\frac{1}{s'}} < \infty ,
\end{gathered}
$\hfill{\rm (2.12)}

{\rm 2)}\hfill
$\displaystyle
\begin{gathered}
b_{k} \le C'  \biggl(\sum\limits_{j=0}^{k}a_{k}\biggr) \quad\text{if}\quad
\sup\limits_{n \in \mathbb{N}}\biggl(\sum\limits_{k=n}^{\infty}\beta^{s}_{k}\biggr)^{\frac{1}{s}}\biggl(\sum\limits_{k=0}^{n}\beta_{k}^{-s'}\biggr)^{\frac{1}{s'}} < \infty.
\end{gathered}
$\hfill{\rm (2.13)}%

\end{theorem}

For $s \in (1,\infty)$ a~proof is given in~\cite{Bennet}. The cases $s=1$ and $s=\infty$ are clear and are verified directly.

\section{The space $B^{\varphi_{0}}_{p,q}(\mathbb{R}^{n},\{t_{k}\})$}

In this section we show that the norms \eqref{Besov1} are equivalent with various functions $\varphi_{0}$ under fairly weak constraints on a~variable smoothness
$\{t_{k}\}$.

 To this aim we shall modify and substantially simplify the methods of~\cite{Ry}. In particular, the use of the theorem on the boundedness of the local maximal Hardy--Littlewood operator
 in weighted Lebesgue spaces (see~\cite{Ry})
 will be avoided. The crux of such a~simplification  in the `correct definition' of the weight class  $X^{\alpha_{3}}_{\alpha,\sigma,p}$ (Definition 2.1) and the
 use of the maximal function $M_{A}(m,j,c)$ (see below).

Let $A > 0$ and let a function $\varphi_{0} \in D(\mathbb{R}^{n})$, $\int\varphi_{0}(x)\,dx=1$, $\varphi=\varphi_{0}-2^{-n}\varphi_{0}(\frac{\cdot}{2})$.
 Given a~distribution $f \in S'_{e}$, consider the following maximal function
\begin{equation}
\label{3.1}
M_{A}(m,j,c)[f]:=\sup\limits_{k \geq j}2^{A(j-k)} \sup\limits_{y \in cQ^{n}_{j,m}} \varphi_{k} \ast f (y), \qquad m \in \mathbb{Z}^{n}, \ \ j \in \mathbb{N}_{0}, \ \ c \geq 1.
\end{equation}

A similar construction was used in~\cite{Ry}. The advantage of our construction is that it is local and proves more natural when working with spaces
of variable smoothness.

Let $c \geq 1$, $m \in \mathbb{Z}^{n}$, $j \in \mathbb{N}_{0}$. Then, as is easily seen, $M_{A}(m,j,c) < \infty$ with  $f \in S'_{e}$, $A \geq N_{f}$.

We shall require the following crucial fact from~\cite{Ry}.

\begin{theorem} [the local reproducing formula]
\label{Th3.1}
Let $\varphi_{0} \in D(\mathbb{R}^{n})$ and $\int\varphi_{0}\,dx = 1$. We set $\varphi:=\varphi_{0}-2^{-n}\varphi_{0}(\frac{\cdot}{2})$. Then, for any number $A \geq 0$,
there exist functions $\psi_{0},\psi \in D(\mathbb{R}^{n})$ such that $L_{\psi} \geq A$ and, for any distribution $f \in D'(\mathbb{R}^{n})$,
\begin{equation}
\label{3.2}
f=\sum\limits_{j=0}^{\infty}\psi_{j} \ast \varphi_{j} \ast f \text{ in the sense of convergence in } D'(\mathbb{R}^{n}).
\end{equation}
\end{theorem}

The estimates from the following two lemmas will be crucial in the subsequent analysis.

\begin{Lm}
\label{Lm3.1}
Let $r \in (0,\infty)$, $c_{1} \geq 1$, $A > 0$. Then there exist constants
$C:=C(n,r,\varphi_{0}, A, c_{1}) > 0$, $c_{2}(\varphi_{0},c_{1},n) \geq c_{1}$ such that, for $j \in \mathbb{N}_{0}$, $m \in \mathbb{Z}^{n}$ and $f \in S'_{e}$,
\begin{equation}
\label{3.3}
 M_{A}(m,j,c_{1})[f] \le C \biggl(\sum\limits_{k=j}^{\infty}2^{(j-k)Ar}2^{kn}\int\limits_{c_{2}Q^{n}_{j,m}}|\varphi_{k}\ast f(z)|^{r}\,dz\biggr)^{\frac{1}{r}}.
\end{equation}
\end{Lm}

The proof of this lemma mainly follows that of Lemma 2.9 of~\cite{Ry}. We give a~sketch of the proof indicating only the differences appearing in our setting.

In view of \eqref{3.2}
$$
\varphi_{j}\ast f = (\psi_{0})_{j} \ast (\varphi_{0})_{j} \ast \varphi_{j}\ast f + \sum\limits_{k=j+1}^{\infty} \varphi_{j} \ast \psi_{k} \ast \varphi_{k} \ast f
$$
(as in the proof of Lemma 2.9 of~\cite{Ry}), where one may assume that $L_{\psi} > A$.

For all $k > j$ we have the following estimate
$$
\|\varphi_{j} \ast \psi_{k} | L_{\infty}(\mathbb{R}^{n})\| \le C 2^{(j-k)A+jn}
$$
with constant $C > 0$ independent of both~$k$ and~$j$
(the condition $L_{\psi} > A$ being useful; see \cite{Ry} for details).

Hence, taking into account that for $k \geq j$ the support of the  function  $\varphi_{j} \ast \psi_{k}$  lies in the cube of side length at most
$\widetilde{c} 2^{-jn}$, we have, for $i \le j$,
\begin{equation}
\label{3.4}
\sup\limits_{y \in c_{1}Q^{n}_{i,m}} \varphi_{j}\ast f(y) \le C \biggl(\sum\limits_{k=j}^{\infty}2^{(j-k)A}2^{kn}\int\limits_{c_{2}Q^{n}_{i,m}}|\varphi_{k}\ast f(z)|\,dz\biggr).
\end{equation}

If $r \geq 1$, then an application of H\"older's inequality first for integrals and then for series with exponents $r,r'$ completes the proof.

If $r \in (0,1)$, then from \eqref{3.4} it clearly follows that, for $k \geq j$,
\begin{equation}
\label{3.5}
\begin{gathered}
2^{(j-k)A}\sup\limits_{y \in c_{1}Q^{n}_{j,m}} \varphi_{k}\ast f(y) \le C \biggl(\sum\limits_{l=k}^{\infty}2^{(j-l)A}2^{ln}\int\limits_{c_{2}Q^{n}_{j,m}}|\varphi_{l}\ast f(z)|\,dz\biggr) \le \\
\le C (M_{A}(m,j,c_{1}))^{1-r}\biggl(\sum\limits_{l=k}^{\infty}2^{(j-l)Ar}2^{ln}\int\limits_{c_{2}Q^{n}_{j,m}}|\varphi_{l}\ast f(z)|^{r}\,dz\biggr).
\end{gathered}
\end{equation}


Now the conclusion of the lemma with $A \geq N_{f}$ follows from \eqref{3.5}, because we have $M_{A}(m,j,c_{1}) < \infty$ for $A \geq N_{f}$.
Note that the constant~$C$ in~\eqref{3.5} depends only on $n,r,\varphi_{0}$, $A$,~$c_{1}$. The proof is completed by arguing, with small modification, as at the end of the proof
of Lemma~2.9 of~\cite{Ry}.

\textbf{Notation}. We need the following notation. Given $\mu \in (0,q]$, we set $q_{\mu}:=\frac{q}{\mu}$ and denote by $q'_{\mu}$
the conjugate H\"older exponent to~$q_{\mu}$. For $r \in (0,p]$  the exponents $p_{r}$ and $p'_{r}$ are defined similarly.

\begin{Lm}
\label{Lm3.3}
Let $p,q \in (0,\infty]$, $r \in (0,p]$,  $\sigma_{1}=r p'_{r}$, $\sigma_{2} \in \mathbb{R}$ and let $\{\alpha^{1}_{k}\}$, $\{\alpha^{2}_{k}\}$ be
sequences of positive real numbers, $\alpha_{3} \geq 0$. Next, let $\{t_{k}\}$ be a~$p$-admissible  weight sequence from $X^{\alpha_{3}}_{\alpha,\sigma,p}$, $c \geq 1$.
Further, assume that, for some numbers $A > 0$, $\mu \in (0,\min\{1,q,r\}]$, condition {\rm (2.12)} is satisfied with $s=q_{\mu}$, $\{\beta_{k}\}=\{(2^{kA}\alpha^{1}_{k})^{\mu}\}$.
Then, for any distribution $f \in B^{\varphi_{0}}_{p,q}(\mathbb{R}^{n},\{t_{k}\})$,
\begin{equation}
\begin{split}
\label{3.6}
   \biggl(\sum\limits_{j=0}^{\infty}\biggl(\sum\limits_{m \in \mathbb{Z}^{n}}t^{p}_{j,m} (M_{A}(m,j,c)[f])^{p}\biggr)^{\frac{q}{p}}\biggr)^{\frac{1}{q}} \le C\biggl(\sum\limits_{k=0}^{\infty}\biggl(\int\limits_{\mathbb{R}^{n}}t^{p}_{k}(x)(\varphi_{k} \ast f(x))^{p}\,dx \biggr)^{\frac{q}{p}}\biggr)^{\frac{1}{q}},
\end{split}
\end{equation}
where the constant $C:=C(n,p,q,\mu,\alpha_{3},\{\alpha^{1}_{k}\},\sigma_{1},c,A,\varphi_{0}) > 0$ is independent of~$f$.
{\rm (}The modifications for $p=\infty$ or $q=\infty$ are clear.{\rm )}
\end{Lm}

The proof will be given for $p,q \neq \infty$ (the case $p=\infty$ or $q=\infty$ is treated similarly).
Using first estimate \eqref{3.3}, and then taking into account that the $l_{q}$-norm is monotone in~$q$,
employing Minkowski's inequality for sums (inasmuch as $\frac{p}{\mu} \geq 1$), and using H\"older's inequality with exponents $p_{r}$ and $p'_{r}$
for integrals, this gives, for  $j \in \mathbb{N}_{0}$ and $\mu \le \min\{1,q,r\}$,
\begin{gather}
\allowdisplaybreaks
\biggl(\sum\limits_{m \in \mathbb{Z}^{n}}t^{p}_{j,m} [M_{A}(m,j,c)]^{p}\biggr)^{\frac{\mu}{p}} \le  \notag \\
\le
C \biggl(\sum\limits_{m \in \mathbb{Z}^{n}}t^{p}_{j,m} \biggl(\sum\limits_{k=j}^{\infty}2^{(j-k)Ar}2^{kn}\int\limits_{\widetilde{c} Q^{n}_{j,m}}|\varphi_{k}\ast f(z)|^{r}\,dz\biggr)^{\frac{\mu p}{\mu r}}\biggr)^{\frac{ r}{p}} \le
\notag\\
\le C \biggl(\sum\limits_{m \in \mathbb{Z}^{n}}t^{p}_{j,m} \biggl(\sum\limits_{k=j}^{\infty}\biggl(2^{(j-k)Ar+kn}\int\limits_{\widetilde{c} Q^{n}_{j,m}}|\varphi_{k}\ast f(z)|^{r}\,dz\biggr)^{\frac{\mu}{r}}\biggr)^{\frac{p}{\mu}}\biggr)^{\frac{\mu}{p}} \le
\notag\\
\le C \sum\limits_{k=j}^{\infty}2^{(j-k)A\mu}\biggl(2^{\frac{knp}{r}}\sum\limits_{m \in \mathbb{Z}^{n}}t^{p}_{j,m}\biggl(\int\limits_{\widetilde{c} Q^{n}_{j,m}}\frac{t^{r}_{k}(x)}{t^{r}_{k}(x)}|\varphi_{k}\ast f(z)|^{r}\,dz\biggr)^{\frac{p}{r}} \biggr)^{\frac{\mu}{p}} \le
\notag\\
\label{3.7}
\le C \sum\limits_{k=j}^{\infty}2^{(j-k)A \mu}\frac{(\alpha^{1}_{j})^{\mu}}{(\alpha^{1}_{k})^{\mu}}\biggl(\sum\limits_{m \in \mathbb{Z}^{n}}\int\limits_{\widetilde{c}Q^{n}_{j,m}}t_{k}^{p}(z)|\varphi_{k}\ast f(z)|^{p}\,dz\biggr)^{\frac{\mu}{p}}.
\end{gather}
Here, we also used that $2^{\frac{knp}{r}}=2^{kn}+2^{kn\frac{p}{\sigma}}$.
Now the required assertion follows from \eqref{3.8} and Theorem \ref{Th2.1} in view of the assumptions on the parameter~$A$ and the sequence $\{\alpha^{1}_{k}\}$.

\begin{theorem}
\label{Th3.2}
Let $p,q \in (0,\infty]$, $r \in (0,p]$, $\sigma_{1}=r p'_{r}$, $\sigma_{2}=p$,  let $\{\alpha^{1}_{k}\}$, $\{\alpha^{2}_{k}\}$ be
sequences of positive real numbers, and let  $\{t_{k}\} \in X^{\alpha_{3}}_{\alpha,\sigma,p}$. Assume that, for some numbers $A > 0$, $\mu \in (0,\min\{1,q,r\}]$ and a~function $\zeta:=\zeta_{0}-2^{-n}\zeta(\frac{\cdot}{2})$  $(\zeta_{0} \in D(\mathbb{R}^{n}))$,  condition {\rm (2.12)} is satisfied with $s=q_{\mu}$, $\{\beta_{k}\}=\{(2^{kA}\alpha^{1}_{k})^{\mu}\}$,
and condition {\rm (2.13)} is satisfied with $s=q_{\mu}$, $\{\beta_{k}\}=\{(2^{-k(1+L_{\zeta})}\alpha^{2}_{k})^{\mu}\}$. Then
\begin{equation}
\label{3.8}
\|f|B^{\zeta_{0}}_{p,q}(\mathbb{R}^{n}, \{t_{k}\})\| \le C \|f|B^{\varphi_{0}}_{p,q}(\mathbb{R}^{n}, \{t_{k}\})\|,
\end{equation}
the constant $C > 0$ is independent of the distribution $f\in S'_e$.
\end{theorem}

The proof will be carried out in the cases $p,q \neq \infty$ (the cases $p=\infty$ and $q=\infty$ are simpler and
are dealt with similarly).

By the local reproducing formula, \begin{equation}
\label{3.9}
f=\sum\limits_{k=0}^{\infty}\psi_{k} \ast \varphi_{k} \ast f  \ \ \text{ in the sense of}\ \ D'(\mathbb{R}^{n}),
\end{equation}
where $L_{\psi}$ may be arbitrarily large.

We set $\zeta:=\zeta-2^{-n}\zeta(\frac{\cdot}{2})$. According to~\cite{Ry},
\begin{equation}
\label{3.10}
\|\zeta_{j} \ast \psi_{k}|L_{\infty}(\mathbb{R}^{n})\| \le C
\begin{cases}
 2^{(j-k)(L_{\psi}+1)}2^{jn}, & j < k \\
 2^{(k-j)(L_{\zeta}+1)}2^{kn}, & k \le j.
\end{cases}
\end{equation}

The diameter of the support of the function $\zeta_{j} \ast \psi_{k}$ is at most $c 2^{-jn}$ (in the case $k > j$) and $c 2^{-kn}$ (in the case $j \geq k$),
the constant~$c$ depending only on $\varphi_{0}$, $\zeta_{0}$, and~$n$. Combining this with \eqref{3.9}, \eqref{3.10}, this establishes
\begin{multline}
\sup\limits_{y \in Q^{n}_{j,m}}|\zeta_{j} \ast f(y)| \le \\
\le \sum\limits_{k=0}^{j}2^{(k-j)(L_{\zeta}+1)}\sup\limits_{y \in cQ^{n}_{k,\widetilde{m}}}|\varphi_{k} \ast f(y)|+\sum\limits_{k=j+1}^{\infty}2^{(j-k)(L_{\psi}+1)}\sup\limits_{y \in cQ^{n}_{j,m}}|\varphi_{k} \ast f(y)|.
\label{3.11}
\end{multline}

For $k < j$, the cube $Q^{n}_{k,\widetilde{m}}$  is the unique dyadic cube  on the right of \eqref{3.11} that has side length $2^{-k}$ and contains the cube $Q^{n}_{j,m}$

Choosing $L_{\psi}$ so that $A < L_{\psi}+1$, we obtain, for $\varepsilon=L_{\psi}+1 - A > 0$,
\begin{equation}
\label{3.12}
\sum\limits_{k=j+1}^{\infty}2^{(j-k)(L_{\psi}+1)}\sup\limits_{y \in cQ^{n}_{j,m}}|\varphi_{k} \ast f(y)| \le \sum\limits_{k=j+1}^{\infty}2^{(j-k)\varepsilon}M_{A}(m,j,c) \le C M_{A}(m,j,c).
\end{equation}

Using  \eqref{3.12} with $\mu \in (0,\min\{1,p,q\}]$ it is found that
\begin{gather}
\biggl(\sum\limits_{m \in \mathbb{Z}^{n}}t^{p}_{j,m}\sup\limits_{y \in Q^{n}_{j,m}}|\zeta_{j} \ast f(y)|^{p}\biggr)^{\frac{1}{p}} \le \notag \\
\le C \biggl(\sum\limits_{m \in \mathbb{Z}^{n}}t^{p}_{j,m}\biggl(\sum\limits_{k=0}^{j}2^{\mu(k-j)(L_{\zeta}+1)}\sup\limits_{y \in cQ^{n}_{k,\widetilde{m}}}|\varphi_{k} \ast f(y)|^{\mu}\biggr)^{\frac{p}{\mu}}\biggr)^{\frac{1}{p}}+
\notag \\
+C \biggl(\sum\limits_{m \in \mathbb{Z}^{n}}t^{p}_{j,m}[M_{A}(m,j,c)]^{p}\biggr)^{\frac{1}{p}}=:S_{1,j}+S_{2,j}.
\label{3.13}
\end{gather}

To estimate $S_{1,j}$ we shall employ Minkowski's inequality for sums (because  $\frac{p}{\mu} \geq 1$) and use condition~(2.3). We have
\begin{gather}
\allowdisplaybreaks
S_{1,j} \le C \biggl(\sum\limits_{m \in \mathbb{Z}^{n}}t^{p}_{j,m}\biggl(\sum\limits_{k=0}^{j}2^{\mu(k-j)(L_{\zeta}+1)}\sup\limits_{y \in cQ^{n}_{k,\widetilde{m}}}|\varphi_{k} \ast f(y)|^{\mu}\biggr)^{\frac{p}{\mu}}\biggr)^{\frac{1}{p}} \le
\notag \\
\le C \biggl(\sum\limits_{k=0}^{j}2^{\mu(k-j)(L_{\zeta}+1)}\biggl(\biggl(\sum\limits_{\substack{\widetilde{m} \in \mathbb{Z}^{n}\\
Q^{n}_{j,m} \subset Q^{n}_{k,\widetilde{m}}}} t^{p}_{j,m}\biggr)\sup\limits_{y \in cQ^{n}_{k,\widetilde{m}}}|\varphi_{k} \ast f(y)|^{p}\biggr)^{\frac{\mu}{p}}\biggr)^{\frac{1}{\mu}} \le
\notag \\
\le C \biggl(\sum\limits_{k=0}^{j}2^{\mu(j-k)(-L_{\zeta}-1)}\biggl(\frac{\alpha^{2}_{j}}{\alpha^{2}_{k}}\biggr)^{\mu}\biggl(\sum\limits_{\widetilde{m} \in \mathbb{Z}^{n}}t^{p}_{k,\widetilde{m}}\sup\limits_{y \in cQ^{n}_{k,\widetilde{m}}}|\varphi_{k} \ast f(y)|^{p}\biggr)^{\frac{\mu}{p}}\biggr)^{\frac{1}{\mu}} \le
\notag \\
\label{3.14}
\le C \biggl(\sum\limits_{k=0}^{j}2^{\mu(j-k)(-L_{\zeta}-1)}\biggl(\frac{\alpha^{2}_{j}}{\alpha^{2}_{k}}\biggr)^{\mu}\biggl(\sum\limits_{m \in \mathbb{Z}^{n}}t^{p}_{k,m}[M_{A}(m,j,c)]^{p}\biggr)^{\frac{\mu}{p}}\biggr)^{\frac{1}{\mu}}.
\end{gather}

Now the conclusion of the theorem follows from estimates \eqref{3.13}, \eqref{3.14}, the restriction on $L_{\zeta}$, and from
Theorem \ref{Th2.1} and Lemma~\ref{Lm3.3}.

\begin{Remark}
\label{R3.2}
{\rm From Theorem \ref{Th3.2} it clearly follows that the spaces $B^{\varphi_{0}}_{p,q}(\mathbb{R}^{n}, \{t_{k}\})$ and $B^{\zeta_{0}}_{p,q}(\mathbb{R}^{n}, \{t_{k}\})$ are equal
and that the corresponding norms are equivalent, provided that (2.13) is satisfied with $s=q_{\mu}$, $\{(2^{-k(1+L_{\zeta})}\alpha^{2}_{k})^{\mu}\}=\{\beta_{k}\}$, and in addition,
with $s=q_{\mu}$, $\{(2^{-k(1+L_{\varphi})}\alpha^{2}_{k})^{\mu}\}=\{\beta_{k}\}$.\marginpar{<--ploho}

In the particular case when $\{\alpha^{1}_{k}\}=\{2^{k\alpha_{1}}\}$, $\{\alpha^{2}_{k}\}=\{2^{k\alpha_{2}}\}$ ($\alpha_{2} \in \mathbb{R}$)
the space $B^{\varphi_{0}}_{p,q}(\mathbb{R}^{n}, \{t_{k}\})$ is independent of the choice of a~function $\varphi_{0} \in D(\mathbb{R}^{n})$,
 provided that $L_{\varphi} \ge [\alpha_{2}]$, the square bracket denoting the integer part.

In the case $p,q \in (0,\infty]$, $p \neq \infty$, $r \in (0,p)$, $\{t_{k}\}=\{2^{ks}\gamma\}$ with $s \in \mathbb{R}$, $\gamma^{p} \in A^{\text{\rm loc}}_{\frac{p}{r}}(\mathbb{R}^{n})$
we obtain one result of~\cite{Ry} (Corollary 2.7)---the independence of the definition of weighted Besov spaces of the choice of a~function $\varphi_{0}$,
provided that $L_{\varphi} \ge  [s]$.}
\end{Remark}

A closer look at the argument given in the proof of Theorem \ref{Th3.2} leads to the following interesting observation.

\begin{theorem}
\label{Th3.3}
Let $p,q \in (0,\infty]$, $r \in (0,p]$, $\sigma_{1}=r p'_{r}$, $\sigma_{2}=p$, let
$\{\alpha^{1}_{k}\}$, $\{\alpha^{2}_{k}\}$ be sequences of positive real numbers, and let $\{t_{k}\} \in X^{\alpha_{3}}_{\alpha,\sigma,p}$.
Assume that, for some numbers $A > 0$, $\mu \in (0,\min\{1,q,r\}]$, condition {\rm (2.12)} is satisfied with $s=q_{\mu}$, $\{\beta_{k}\}=\{(2^{kA}\alpha^{1}_{k})^{\mu}\}$,
and condition {\rm (2.13)} is satisfied with $s=q_{\mu}$, $\{\beta_{k}\}=\{(2^{-k(1+L_{\varphi})}\alpha^{2}_{k})^{\mu}\}$.

Then $B^{\varphi_{0}}_{p,q}(\mathbb{R}^{n},\{t_{k}\})=B^{\varphi_{0}}_{p,q}(\mathbb{R}^{n},\{\overline{t}_{k}\})$, the corresponding norms being equivalent.
\end{theorem}

The proof requires only minor modifications to that of Theorem~\ref{Th3.2}, the arguments depend upon Lemma \ref{Lm2.1}. We omit the details.

\begin{Remark} \rm
Theorem \ref{Th3.3} enables us to get rid of the singularities and degeneracy points in weights from a sequence $\{t_{k}\}$.
This fact proves to be a~useful tool in subsequent applications of the spaces $B^{\varphi_{0}}_{p,q}(\mathbb{R}^{n},\{t_{k}\})$.
In particular, if $p,q \in (0,\infty]$, $p \neq \infty$, $\gamma^{p} \in A^{\text{\rm loc}}_{\infty}(\mathbb{R}^{n})$, $L_{\varphi} \ge [s]$, then $B^{s}_{p,q}(\mathbb{R}^{n},\gamma)=B^{\varphi_{0}}_{p,q}(\mathbb{R}^{n},\{\gamma_{k}\})$ where $\gamma_{k}(x)=2^{ks+\frac{kn}{p}}\sum\limits_{m \in \mathbb{Z}^{n}}\chi_{\widetilde{Q}^{n}_{k,m}}(x)\|\gamma|L_{p}(\widetilde{Q}^{n}_{k,m})\|$ with $k \in \mathbb{N}_{0}$, $x \in \mathbb{R}^{n}$.
\end{Remark}

\section{Comparison of the spaces $\widetilde{B}^{l}_{p,q,r}(\mathbb{R}^{n},\{t_{k}\})$ and $B^{\varphi_{0}}_{p,q}(\mathbb{R}^{n},\{t_{k}\})$.}

It is well known (see, for example, \S\,2.5.3 of~\cite{Tri}) that for $\{s_{k}\}=\{2^{ks}\}$ and $p,q \in (0,\infty]$
the spaces $\widetilde{B}^{l}_{p,q,p}(\mathbb{R}^{n},\{s_{k}\})$ and $B^{\{s_k\}}_{p,q}(\mathbb{R}^{n})$ coincide when
$s > n\max\{\frac{1}{p}-1,0\}$. However, if $s < n\max\{\frac{1}{p}-1,0\}$, then the space $B^{\{s_k\}}_{p,q}(\mathbb{R}^{n})$
(unlike the space  $\widetilde{B}^{l}_{p,q,p}(\mathbb{R}^{n},\{t_{k}\})$) contains
the Dirac delta function. In this section we shall establish the relation between the spaces $\widetilde{B}^{l}_{p,q,r}(\mathbb{R}^{n},\{t_{k}\})$ and
$B^{\varphi_{0}}_{p,q}(\mathbb{R}^{n},\{t_{k}\})$ under minimal assumptions
on the variable smoothness $\{t_{k}\}$.

\begin{theorem}
\label{Th4.1}
Let $p,q \in (0,\infty]$, $r \in [1,\infty]$, $\theta \in (0,\min\{p,r\}]$, $\sigma_{1}=\theta p'_{\theta}$, $\sigma_{2} \in (0,\infty]$, $\alpha_{3} \geq 0$.
Assume that the sequence $\{2^{kn\mu (\frac{1}{\theta}-\frac{1}{r}) }(\alpha^{1}_{k})^{-\mu}\}$ lies in $l_{q'_{\mu}}$ for some
$\mu \le \min\{1,q,\theta\}$. Further assume that a~weight sequence $\{t_{k}\} \in X^{\alpha_{3}}_{\alpha,\sigma,p}$. Then
$B^{\varphi_{0}}_{p,q}(\mathbb{R}^{n},\{t_{k}\}) \subset L^{\text{\rm loc}}_{r}(\mathbb{R}^{n})$.
\end{theorem}

The proof will be conducted for  $p,q \neq \infty$. The cases $p=\infty$ or $q=\infty$ are dealt with similarly.

Let $f \in B^{\varphi_{0}}_{p,q}(\mathbb{R}^{n},\{t_{k}\})$. From Minkowski's inequality and since the $l_{q}$-norm is monotone in~$q$, we have, for $\mu \le \min\{1,q,\theta\}$, $m\in \mathbb Z^n$,
\begin{equation}
\begin{gathered}
\label{4.1}
I_{m}:=\sum\limits_{j=0}^{\infty}\biggl(\int\limits_{Q^{n}_{0,m}}|\varphi_{j}\ast f(y)|^{r}\,dy\biggr)^{\frac{1}{r}} \le
\biggl[\sum\limits_{j=0}^{\infty}\biggl(\int\limits_{Q^{n}_{0,m}}|\varphi_{j}\ast f(y)|^{r}\,dy\biggr)^{\frac{\mu}{r}}\biggr]^{\frac{1}{\mu}} \le \\
\le \biggl[\sum\limits_{j=0}^{\infty}\biggl(\sum\limits_{\substack{\widetilde{m} \in \mathbb{Z}^{n}\\
Q^{n}_{j,\widetilde{m}} \subset Q^{n}_{0,m}}}\biggl(\int\limits_{Q^{n}_{j,\widetilde{m}}}|\varphi_{j}\ast f(y)|^{r}\,dy\biggr)^{\frac{\theta}{r}}\biggr)^{\frac{\mu}{\theta}}\biggr]^{\frac{1}{\mu}}
\end{gathered}
\end{equation}

Using Lemma \ref{Lm3.1} (with $\theta$ instead of~$r$),  condition (2.3) and  H\"older's inequality for integrals, this establishes
\begin{gather}
\allowdisplaybreaks
\sum\limits_{\substack{\widetilde{m} \in \mathbb{Z}^{n}\\
Q^{n}_{j,\widetilde{m}} \subset Q^{n}_{0,m}}}\biggl(\int\limits_{Q^{n}_{j,\widetilde{m}}}|\varphi_{j}\ast f(y)|^{r}\,dy\biggr)^{\frac{\theta}{r}}
\le 2^{-nj\frac{\theta}{r}}\sum\limits_{\substack{\widetilde{m} \in \mathbb{Z}^{n}\\
Q^{n}_{j,\widetilde{m}} \subset Q^{n}_{0,m}}}\sup\limits_{y \in Q^{n}_{j,\widetilde{m}}}|\varphi_{j}\ast f(y)|^{\theta} \le
\notag\\
\le C 2^{j\theta(A-\frac{n}{r})}\sum\limits_{\substack{\widetilde{m} \in \mathbb{Z}^{n}\\
Q^{n}_{j,\widetilde{m}} \subset Q^{n}_{0,m}}}\sum\limits_{k=j}^{\infty}2^{k\theta(\frac{n}{\theta}-A)}\int\limits_{cQ^{n}_{j,\widetilde{m}}}|\varphi_{k}\ast f(z)|^{\theta}\,dz \le
\notag\\
\le C 2^{j\theta(A-\frac{n}{r})} \sum\limits_{k=j}^{\infty}2^{k\theta(\frac{n}{\theta}-A)}\biggl[\int\limits_{cQ^{n}_{0,m}}t^{p}_{0}(z)\,dz\biggr]^{\frac{\theta}{p}}\times \notag \\
\times \biggl[\int\limits_{c Q^{n}_{0,m}}(t_{k}(z))^{-\theta p'_{\theta}}\biggr]^{\frac{1}{p'_{\theta}}}\biggl[\int\limits_{cQ^{n}_{0,m}}t^{p}_{k}(z)|\varphi_{k}\ast f(z)|^{p}\,dz\biggr]^{\frac{\theta}{p}} \le
\notag\\
\label{4.2}
\le C 2^{j\theta(A-\frac{n}{r})} \sum\limits_{k=j}^{\infty}2^{k\theta(\frac{n}{\theta}-A)}\biggl(\frac{\alpha^{1}_{0}}{\alpha^{1}_{k}}\biggr)^{\theta}\biggl[\int\limits_{cQ^{n}_{0,m}}t^{p}_{k}(z)|\varphi_{k}\ast f(z)|^{p}\,dz\biggr]^{\frac{\theta}{p}}.
\end{gather}
The constant $C$ on the right  of \eqref {4.2} depends on~$m\in \mathbb Z^n$.

We first substitute estimate \eqref{4.2} into \eqref{4.1}. Then, for $A > \frac{n}{\theta}$,  we change the order of summation, apply H\"older's
inequality with exponents $q_{\mu}$, $q'_{\mu}$ for sums
(here, we set $q_{\mu}=\frac{q}{\mu}$), and use the condition $\{2^{kn\mu(\frac{1}{\theta}-\frac{1}{r})}(\alpha^{1}_{k})^{-\mu}\} \in l_{q'_{\mu}}$. Hence,
\begin{equation}
\label{4.3}
\begin{gathered}
I^{q}_{m} \le C \biggl(\sum\limits_{j=0}^{\infty}\biggl(\frac{2^{kn(\frac{1}{\theta}-\frac{1}{r})}}{\alpha^{1}_{k}}\biggr)^{\mu}\biggl[\int\limits_{cQ^{n}_{0,m}}t^{p}_{k}(z)|\varphi_{k}\ast f(z)|^{p}\,dz\biggr]^{\frac{\mu}{p}}\biggr)^{\frac{q}{\mu}} \le \\
 \le C \sum\limits_{k=0}^{\infty}\biggl(\int\limits_{cQ^{n}_{0,\widetilde{m}}}t^{p}_{k}(z)|\varphi_{k}\ast f(z)|^{p}\,dz\biggr)^{\frac{q}{p}}.
\end{gathered}
\end{equation}

From \eqref{4.3} it follows that the series $\sum\limits_{j=1}^{\infty}\varphi_{j}\ast f$ converges in the sense of
 $L^{\text{\rm loc}}_{r}(\mathbb{R}^{n})$ to some function $g \in L^{\text{\rm loc}}_{r}(\mathbb{R}^{n})$.
 Hence, since $r \geq 1$, this shows that the series $\sum\limits_{j=1}^{\infty}\varphi_{j}\ast f$ converges to~$g$ in $L^{\text{\rm loc}}_{1}(\mathbb{R}^{n})$,
 and hence in the sense of $D'(\mathbb{R}^{n})$.
But then $f=g \in L^{\text{\rm loc}}_{r}(\mathbb{R}^{n})$, because, clearly, the series $\sum\limits_{j=1}^{\infty}\varphi_{j}\ast f$
converges to~$f$ in the sense of $D'(\mathbb{R}^{n})$. The proof is complete.

For further purposes we shall need to recall some elements of spline approximation theory.
By $\Sigma^{l}_{k}$ we shall denote the linear space of splines of degree~$l$ of the form
$$
S(x):=\sum\limits_{m \in \mathbb{Z}^{n}} \beta_{k,m} N^{l}_{k,m}(x) , \qquad x\in \mathbb{R}^{n},
$$
where $N^{l}_{k,m}$ is the dyadic $B$-spline of degree~$l$.
For the background definitions and references we refer the reader to~\cite{Tyu}.

We shall need the following basic properties of the $B$-splines  $N^{l}_{k,m}$:

1) The $B$-splines $N^{l}_{k,m}$ form a~partition of unity on~$\mathbb{R}^{n}$ for each fixed
$k \in \mathbb{N}_{0}$. That is,
$$
\sum\limits_{m \in \mathbb{Z}^{n}}N^{l}_{k,m}(x)=1, \qquad x \in \mathbb{R}^{n}.
$$
Here, the overlapping multiplicity of the supports of splines $N^{l}_{k,m}$ is finite
and is independent of both~$k$ and~$m$. We also note that $\operatorname{supp} N^{l}_{k,m} \subset
\frac{m}{2^{k}}+[0,\frac{l}{2^{k}}]^{n}$ and $N^{l-1}_{k,m}(x) \in (0,1]$ for $x
\in \frac{m}{2^{k}}+(0,\frac{l}{2^{k}})^{n}$.

2) On each cube $Q^{n}_{k,m}$ the function $N^{l}_{k,m}$ is a polynomial of degree $\le l$ in each variable.

3) The spline $N^{l}$ has continuous derivative of order $l-1$. At knots $t_{i}=i$, $i \in \{0,1,\dots,l\}$, the spline $N^{l}$
has finite one-sided derivatives of order $l$. Hence, for some $C>0$ (which is independent of $x,h,k$),
$$
\Delta^{l}(h)N^{l}_{k,m}(x) \le C (2^{k}|h|)^{l}, \qquad x,h \in \mathbb{R}^{n}.
$$

4) Any spline $S=\sum\limits_{m \in Z^{n}}\beta_{k,m}N^{l}_{k,m}$ may be expanded in a~series of splines $N^{l}_{j,m}$ for $j \geq k$:
$$
  S=\sum\limits_{m \in Z^{n}}\widehat{\beta}_{k,m}(S)N^{l}_{j,m}.
$$

\begin{Lm}
\label{Lm3.2}
Let $l\in \mathbb{N}$ and let a spline $S \in \Sigma^{l}_{k}$. Then, for any $r_{1}, r_{2} \in (0,+\infty]$ and any cube $Q^{n}_{k,m}$,
\begin{multline}
C_{1}\|S|L_{r_{1}}(Q^{n}_{k,m})\| \le
\biggl(\!\!\!\sum\limits_{\substack{\widetilde{m} \in \mathbb{Z}^{n}\\
Q^{n}_{k,m} \bigcap \operatorname{supp} N^{l-1}_{k,\widetilde{m}} \neq
\emptyset}}|\alpha_{k,\widetilde{m}}(S)|^{r_{1}}2^{-kn}\biggr)^{\frac{1}{r_{1}}} \le \\
\le C_{2}2^{kn(\frac{1}{r_{2}}-\frac{1}{r_{1}})}\|S|L_{r_{2}}(C_{3}Q^{n}_{k,m})\|,
\label{eq4.4} 
\end{multline}
the constants $C_{1},C_{2},C_{3} > 0$ are independent of both the cube $Q^{n}_{k,m}$ and the spline~$S$.
\end{Lm}

The proof is a straightforward modification of that of Lemma 4.2 in \cite{DeVore}.

\begin{theorem}
\label{Th4.2}
Let $p,q,r \in (0,\infty]$, $\theta \in (0,\min\{p,r\}]$, $\alpha_{3} \geq 0$, $\sigma_{1}=\theta p'_{\theta}$, $\sigma_{2}=p$,
 and let sequences $\{\alpha^{1}_{k}\}$, $\{\alpha^{2}_{k}\}$ of positive numbers be such that, for some $\mu \in (0,\min\{1,\theta,q\}]$,
 condition {\rm (2.12)} is satisfied with $s=q_{\mu}$, $\{\beta_{k}\}=\{(\frac{\alpha^{1}_{k}}{2^{kn(\frac{1}{\theta}-\frac{1}{r})}})^{\mu}\}$,
 and condition {\rm (2.13)} is satisfied with $s=q_{\mu}$, $\{\beta_{k}\}=\{(2^{-kl}\alpha^{2}_{k})^{\mu}\}$. Next,
 let $\{t_{k}\}$  be a~$p$-admissible a~weight sequence such that $\{\overline{t}_{k}\} \in X^{\alpha_{3}}_{\alpha,\sigma,p}$. Then:

{\rm 1)} Each function $f \in \widetilde{B}^{l}_{p,q,r}(\mathbb{R}^{n},\{t_{k}\})$ can be expanded in an $L^{\text{\rm loc}}_{r}(\mathbb{R}^{n})$-convergent  series of splines. More precisely,
\begin{gather}
f=\sum\limits_{k=0}^{\infty}v^{l}_{k}(f) \ \ \text{ in the sense of} \ L_{r}^{\text{\rm loc}}({R}^{n}), \ \text{ where }
\notag \\
\label{3.37}
v^{l}_{k}(\varphi)(x)=\sum\limits_{m \in \mathbb{Z}^{n}}\beta_{k,m}(\varphi)N^{l}_{k,m}(x), \qquad x \in \mathbb{R}^{n}.
\end{gather}

Moreover, {\rm (}with the standard modifications for $p=\infty$ or $q=\infty)$
\begin{equation}
\label{3.38}
N(f,l,\{t_{k}\}):=\inf\biggl(\sum\limits_{k=0}^{\infty}\biggl(\sum\limits_{m \in \mathbb{Z}^{n}}t^{p}_{k,m}|\beta_{k,m}|^{p}\biggr)^{\frac{q}{p}}\biggr)^{\frac{1}{q}} \le  C \|f|\widetilde{B}_{p,q,r}^{l}(\mathbb{R}^{n},\{t_{k}\})\|,
\end{equation}
the infimum on the left of \eqref{3.38} is taken over all series of splines that $L^{\text{\rm loc}}_{r}(\mathbb{R}^{n})$-converge to~$f$.

{\rm 2)} If, for some multiple sequence $\{\beta_{k,m}\}$,
$$
\biggl(\sum\limits_{k=0}^{\infty}\biggl(\sum\limits_{m \in \mathbb{Z}^{n}}t^{p}_{k,m}|\beta_{k,m}|^{p}\biggr)^{\frac{q}{p}}\biggr)^{\frac{1}{q}} < \infty,
$$
then the series $\sum\limits_{k=0}^{\infty}\sum\limits_{m \in \mathbb{Z}^{n}}\beta_{k,m}N^{l}_{k,m}$ converges in $L_{r}^{\text{\rm loc}}(\mathbb{R}^{n})$ to some function
$f \in \widetilde{B}_{p,q,r}^{l}(\mathbb{R}^{n},\{t_{k}\})$, and moreover, for some constant $C > 0$,
$$
  \|f|\widetilde{B}_{p,q,r}^{l}(\mathbb{R}^{n},\{t_{k,m}\})\| \le C N(f,l,\{t_{k}\}).
$$
\end{theorem}

\textbf{Proof.} The proof of Theorem~\ref{Th4.2} is similar to that of Corollary 4.4 of~\cite{Tyu}.

\begin{theorem}
\label{Th4.3}
Let $p,q \in (0,\infty]$, $r \in (0,p]$, $\alpha_{3} \geq 0$, $\sigma_{1}=r p'_{r}$, $\sigma_{2}=p$, $\varphi_{0} \in D(\mathbb{R}^{n})$, $\varphi:=\varphi_{0}-2^{-n}\varphi_{0}(\frac{\cdot}{2})$,
and let sequences $\{\alpha^{1}_{k}\}$, $\{\alpha^{2}_{k}\}$ of positive numbers  be such that, for some $\mu \in (0,\min\{1,r,q\}]$, condition {\rm (2.12)}
is satisfied with $s=q_{\mu}$, $\{\beta_{k}\}=\{(\alpha^{1}_{k})^{\mu}\}$ and condition {\rm (2.13)} is satisfied with $s=q_{\mu}$, $\{\beta_{k}\}=\{(2^{-k(1+L_{\varphi})}\alpha^{2}_{k})^{\mu}\}$ and
with $s=q_{\mu}$, $\{\beta_{k}\}=\{(2^{-kl}\alpha^{2}_{k})^{\mu}\}$. Next, let  $\{t_{k}\} $ be a~$p$-admissible weight sequence, $\{t_{k}\} \in X^{\alpha_{3}}_{\alpha,\sigma,p}$.

Then  $B^{\varphi_{0}}_{p,q}(\mathbb{R}^{n}, \{t_{k}\})\bigcap L^{\text{\rm loc}}_{1}(\mathbb{R}^{n})=\widetilde{B}^{l}_{p,q,r}(\mathbb{R}^{n},\{t_{k}\})\bigcap L^{\text{\rm loc}}_{1}(\mathbb{R}^{n})$,
the corresponding norms being equivalent.
\end{theorem}

\textbf{Proof.} Throughout the proof we fix $\mu \in (0,\min\{1,r,q\}]$ from the hypotheses of the theorem.

\textit{Step 1}. We claim that $B^{\varphi_{0}}_{p,q}(\mathbb{R}^{n}, \{t_{k}\})\bigcap L^{\text{\rm loc}}_{1}(\mathbb{R}^{n}) \subset \widetilde{B}^{l}_{p,q,r}
(\mathbb{R}^{n},\{t_{k}\})\bigcap L^{\text{\rm loc}}_{1}(\mathbb{R}^{n})$. Indeed, since, by the condition $f \in L^{\text{\rm loc}}_{1}(\mathbb{R}^{n})$,
a~function~$f$ can be expanded in a~series that converges not only in~$D'(\mathbb{R}^{n})$, but also in  $L^{\text{\rm loc}}_{1}(\mathbb{R}^{n})$:
$$
f=\sum\limits_{j=0}^{\infty}\psi_{j} \ast \varphi_{j}\ast f.
$$

We set $f_{1,k}:=\sum\limits_{j=0}^{k}\psi_{j} \ast \varphi_{j}\ast f$, $f_{2,k}:=f-f_{1,k}$, where $k \in \mathbb{N}_{0}$.

We have
\begin{gather*}
\biggl(\sum\limits_{k=0}^{\infty}\biggl(\sum\limits_{m \in \mathbb{Z}^{n}}t^{p}_{k,m}[\delta^{l}_{r}(Q^{n}_{k,m})f]^{p}\biggr)^{\frac{q}{p}}\biggr)^{\frac{1}{q}} \le C \biggl(\sum\limits_{k=0}^{\infty}\biggl(\sum\limits_{m \in \mathbb{Z}^{n}}t^{p}_{k,m}[\delta^{l}_{r}(Q^{n}_{k,m})f_{1,k}]^{p}\biggr)^{\frac{q}{p}}\biggr)^{\frac{1}{q}} + \\
\\
+C \biggl(\sum\limits_{k=0}^{\infty}\biggl(\sum\limits_{m \in \mathbb{Z}^{n}}t^{p}_{k,m}[\delta^{l}_{r}(Q^{n}_{k,m})f_{2,k}]^{p}\biggr)^{\frac{q}{p}}\biggr)^{\frac{1}{q}}=:S_{1}+S_{2}.
\end{gather*}

Since the $l_{q}$-norm is monotone in $q$, we have
\begin{equation}
\label{4.5}
S_{2} \le \biggl(\sum\limits_{k=0}^{\infty}\biggl(\sum\limits_{m \in \mathbb{Z}^{n}}t^{p}_{k,m}\biggl\{\sum\limits_{j=k}^{\infty}[\delta^{l}_{r}(Q^{n}_{k,m})\psi_{j} \ast \varphi_{j} \ast f]^{\mu}\biggr\}^{\frac{p}{\mu}}\biggr)^{\frac{q}{p}}\biggr)^{\frac{1}{q}}.
\end{equation}

To estimate the right-hand side of \eqref{4.5} we shall employ Minkowski's inequality for sums (because $\frac{p}{\mu} \geq 1$).
As a result,
\begin{equation}\label{4.6}
\begin{gathered}
S_{2} \le \biggl(\sum\limits_{k=0}^{\infty}\biggl\{\sum\limits_{j=k}^{\infty}\biggl(\sum\limits_{m \in \mathbb{Z}^{n}}t^{p}_{k,m}[\delta^{l}_{r}(Q^{n}_{k,m})\psi_{j}
\ast \varphi_{j} \ast f]^{p}\biggr)^{\frac{\mu}{p}}\biggr\}^{\frac{q}{\mu}}\biggr)^{\frac{1}{q}} \le \\
\le \biggl(\sum\limits_{k=0}^{\infty}\biggl\{\sum\limits_{j=k}^{\infty}\biggl(\sum\limits_{m \in \mathbb{Z}^{n}}t^{p}_{k,m}2^{\frac{knp}{r}}\|\psi_{j} \ast \varphi_{j} \ast f|L_{r}(cQ^{n}_{k,m})\|^{p}\biggr)^{\frac{\mu}{p}}\biggr\}^{\frac{q}{\mu}}\biggr)^{\frac{1}{q}}
\end{gathered}
\end{equation}
with some constant $c>0$ depending only on~$l$,$n$

Using Holder inequality for integrals, (2.3), \eqref{2.5} and since the cubes $c Q^{n}_{k,m}$ have finite overlapping multiplicity  (independent of~$k$, $m$), this establishes
 \begin{equation}
\begin{split}
\label{4.7}
&\sum\limits_{m \in \mathbb{Z}^{n}}t^{p}_{k,m}2^{\frac{knp}{r}}\|\psi_{j} \ast \varphi_{j} \ast f|L_{r}(cQ^{n}_{k,m})\|^{p} \le C \sum\limits_{m \in \mathbb{Z}^{n}} \frac{\alpha^{1}_{k}}{\alpha^{1}_{j}}\int\limits_{cQ^{n}_{k,m}}t^{p}_{j}(x)|\varphi_{j} \ast \psi_{j} \ast f (x)|^{p} \,dx \le \\
&\le C
\frac{\alpha^{1}_{k}}{\alpha^{1}_{j}}\sum\limits_{\substack{\widetilde{m} \in \mathbb{Z}^{n}\\
Q^{n}_{j,\widetilde{m}} \subset \widetilde{c}Q^{n}_{k,m}}}t^{p}_{j,\widetilde{m}}\sup\limits_{y \in  Q^{n}_{j,\widetilde{m}}} |\varphi_{j} \ast f(y)|^{p} \le C \frac{\alpha^{1}_{k}}{\alpha^{1}_{j}} \sum\limits_{m \in \mathbb{Z}^{n}} t^{p}_{j,m}\sup\limits_{y \in  Q^{n}_{j,m}} |\varphi_{j} \ast f(y)|^{p}.
\end{split}
\end{equation}

Next, inserting \eqref{4.7} into \eqref{4.6}, using Theorem \ref{Th2.1}, and then employing Lemma \ref{Lm3.3}, we obtain
\begin{gather}
S_{2} \le  C \biggl(\sum\limits_{j=0}^{\infty}\biggl(\sum\limits_{m \in \mathbb{Z}^{n}}t^{p}_{j,m} \sup\limits_{y \in c Q^{n}_{j,m}} |\varphi_{j} \ast f(y)|^{p}\biggr)^{\frac{q}{p}}\biggr)^{\frac{1}{q}} \le
 \notag \\
 \le C \biggl(\sum\limits_{j=0}^{\infty}\biggl(\sum\limits_{m \in \mathbb{Z}^{n}}t^{p}_{j,m} [M_{A}(m,j,c)]^{p}\biggr)^{\frac{q}{p}}\biggr)^{\frac{1}{q}} \le
C \|f|B^{\varphi_{0}}_{p,q}(\mathbb{R}^{n},\{t_{k}\})\|.
\label{4.8}
\end{gather}

Since the $l_{q}$-norm is monotone in $q$, we have
\begin{equation}
\label{4.9}
S_{1} \le \biggl(\sum\limits_{k=0}^{\infty}\biggl(\sum\limits_{m \in \mathbb{Z}^{n}}t^{p}_{k,m}\biggl\{\sum\limits_{j=0}^{k}[\delta^{l}_{r}(Q^{n}_{k,m})\psi_{j} \ast \varphi_{j} \ast f]^{\mu}\biggr\}^{\frac{p}{\mu}}\biggr)^{\frac{q}{p}}\biggr)^{\frac{1}{q}}.
\end{equation}

Given $j \le k$ we let $Q^{n}_{j,\widetilde{m}}$ denote the unique dyadic cube containing the cube $Q^{n}_{k,m}$.
By the Lagrange mean-value theorem we have, for $x \in Q^{n}_{j,m}$, $h \in 2^{-k}I^{n}$,
\begin{equation}
\begin{gathered}
\label{4.10}
\Delta^{l}(h)[2^{-jn}\psi_{j} \ast \varphi_{j} \ast f](x) \le |h|^{l}\biggl(\sum\limits_{|\alpha|=l}\sup\limits_{\xi \in \mathbb{R}^{n}}2^{-jn} D^{\alpha}\psi_{j}(\xi)\biggr)\sup\limits_{y \in cQ^{n}_{j,m}}|\varphi_{j} \ast f(y)| \le \\
\le C 2^{(j-k)l} \sup\limits_{y \in cQ^{n}_{j,m}}|\varphi_{j} \ast f(y)|.
\end{gathered}
\end{equation}

Now to estimate the right-hand side of \eqref{4.9} it suffices to invoke Minkowski's inequality for sums (this is possible, because $\frac{p}{\mu} \geq 1$),
and then apply estimate \eqref{4.10}, Theorem~\ref{Th2.1}, and finally, Lemma~\ref{Lm3.3}. We have
\begin{gather}
\allowdisplaybreaks
S_{1} \le \biggl(\sum\limits_{k=0}^{\infty}\biggl(\sum\limits_{j=0}^{k}\biggl\{\sum\limits_{m \in \mathbb{Z}^{n}}t^{p}_{k,m}[\delta^{l}_{r}(Q^{n}_{k,m})\psi_{j} \ast \varphi_{j} \ast f]^{p}\biggr\}^{\frac{\mu}{p}}\biggr)^{\frac{q}{\mu}}\biggr)^{\frac{1}{q}} \le
\notag \\
\le C \biggl(\sum\limits_{k=0}^{\infty}\biggl(\sum\limits_{j=0}^{k}\biggl\{\sum\limits_{m \in \mathbb{Z}^{n}}t^{p}_{k,m} 2^{(j-k)lp}\sup\limits_{y \in cQ^{n}_{j,\widetilde{m}}}|\varphi_{j} \ast f(y)|^{p}\biggr\}^{\frac{\mu}{p}}\biggr)^{\frac{q}{\mu}}\biggr)^{\frac{1}{q}} \le
\notag  \\
 \le  C \biggl(\sum\limits_{k=0}^{\infty}\biggl(\biggl(\frac{\alpha^{2}_{k}}{2^{kl}}\biggr)^{\mu} \sum\limits_{j=0}^{k}\biggl(\frac{2^{jl}}{\alpha^{2}_{j}}\biggr)^{\mu}\biggl\{\sum\limits_{\widetilde{m} \in \mathbb{Z}^{n}}t^{p}_{j,m} \sup\limits_{y \in cQ^{n}_{j,\widetilde{m}}}|\varphi_{j} \ast f(y)|^{p}\biggr\}^{\frac{\mu}{p}}\biggr)^{\frac{q}{\mu}}\biggr)^{\frac{1}{q}} \le
 \notag \\
\label{4.11}
   \le C \biggl(\sum\limits_{j=0}^{\infty}\biggl(\sum\limits_{m \in \mathbb{Z}^{n}} \sup\limits_{y \in cQ^{n}_{j,m}}|\varphi_{j} \ast f(y)|^{p}\biggr)^{\frac{q}{p}}\biggr)^{\frac{1}{q}} \le C \|f|B^{\varphi_{0}}_{p,q}(\mathbb{R}^{n}, \{t_{k}\})\|.
\end{gather}

The required embedding now follows from estimates \eqref{4.6} and \eqref{4.11}.

\textit{Step 2}. We claim that $\widetilde{B}^{l}_{p,q,r}(\mathbb{R}^{n},\{t_{k}\})\bigcap L^{\text{\rm loc}}_{1}(\mathbb{R}^{n}) \subset B^{\varphi_{0}}_{p,q}(\mathbb{R}^{n}, \{t_{k}\})\bigcap L^{\text{\rm loc}}_{1}(\mathbb{R}^{n})$.
Let us expand a~function $f \in \widetilde{B}^{l}_{p,q,r}(\mathbb{R}^{n},\{t_{k}\})\bigcap L^{\text{\rm loc}}_{1}(\mathbb{R}^{n})$ in a~series of splines
converging in  $L^{\text{\rm loc}}_{r}(\mathbb{R}^{n})$. More precisely, we expand $f$ as
\begin{equation}
\begin{gathered}\label{4.12}
f=\sum\limits_{j=0}^{\infty} v^{l}_{j}(f) \quad \text{in } \ \ L^{\text{\rm loc}}_{r}(\mathbb{R}^{n}), \quad \text{ where } \\
v^{l}_{j}(x):=\sum\limits_{m \in \mathbb{Z}^{n}} \beta_{j,m}(v^{l}_{j})N^{l}_{j,m}(x), \qquad x \in \mathbb{R}^{n}.
\end{gathered}
\end{equation}

Let $\omega \in C^{\infty}_{0}$, $\int\omega(x)\,dx = 1$. We set $\Omega_{0}(x):=-\sum\limits_{i=1}^{l}(-1)^{l-i}C^{i}_{l}\frac{1}{i^{n}}\omega(\frac{x}{i})$ for
$x \in \mathbb{R}^{n}$. It is easily seen that $\int\Omega_{0}(x)\,dx=1$ (here, the identity $0=(1-1)^l=\sum\limits_{i=0}^{l}(-1)^{l-i}C^{i}_{l}$ is useful).
We set $\Omega:=\Omega_{0}-2^{-n}\Omega_{0}(\frac{\cdot}{2})$. Hence,
\begin{gather}
\Omega_{j} \ast f = \sum\limits_{i=1}^{l}(-1)^{l-i}C^{i}_{l}\int \frac{2^{jn}}{i^{n}} \biggl[\omega\biggl(2^{j}\biggl(\frac{x-y}{i}\biggr)\biggr) - 2^{-n}\omega\biggl(2^{j-1}\biggl(\frac{x-y}{i}\biggr)\biggr)\biggr]f(y)\,dy = \notag \\
 =\sum\limits_{i=1}^{l}(-1)^{l-i}C^{i}_{l}\biggl(\frac{2^{jn}}{i^{n}}\int\omega(2^{j}y)f(x-iy)\,dy-
\frac{2^{(j-1)n}}{i^{n}}\int\omega(2^{j-1}y)f(x-iy)\,dy\biggr)=
\notag \\
=\int2^{jn}\omega(2^{j}y)\Delta^{l}(y)f(x)\,dy - \int2^{(j-1)n}\omega(2^{(j-1)}y)\Delta^{l}(y)f(x)\,dy;
\label{4.13}
\end{gather}
 the convolution can be understood in the conventional sense, because $f \in L^{\text{\rm loc}}_{1}(\mathbb{R}^{n})$.

The subsequent argument is closely similar to that of Step~1, so we shall give a sketch of this proof, omitting the details.

We set $V^{l}_{1,k}:=\sum\limits_{j=0}^{k}v^{l}_j$, $V^{l}_{2,k}:=f-V^{l}_{1,k}$ for $k \in \mathbb{N}_{0}$. By Theorem~\ref{Th3.2},
\begin{gather}
\biggl(\sum\limits_{k=0}^{\infty}\biggl(\int t^{p}_{k}(x)|\varphi_{k} \ast f(y)|^{p}\biggr)^{\frac{q}{p}}\biggr)^{\frac{1}{q}} \le C \biggl(\sum\limits_{k=0}^{\infty}\biggl(\int t^{p}_{k}(x)|\Omega_{k} \ast f(y)|^{p}\biggr)^{\frac{q}{p}}\biggr)^{\frac{1}{q}} \le \notag \\
\le C \biggl(\sum\limits_{k=0}^{\infty}\biggl(\int t^{p}_{k}(x)|\Omega_{k} \ast V^{l}_{1,k}(y)|^{p}\biggr)^{\frac{q}{p}}\biggr)^{\frac{1}{q}} + C \biggl(\sum\limits_{k=0}^{\infty}\biggl(\int t^{p}_{k}(x)|\Omega_{k} \ast V^{l}_{2,k}(y)|^{p}\biggr)^{\frac{q}{p}}\biggr)^{\frac{1}{q}}=: \notag \\
=:S_{1}+S_{2}
\label{4.14}
\end{gather}
(we note that this theorem is valid without any constrains on the moment of the function~$\Omega$).

Using Minkowski's inequality for sums,
\begin{equation}
\label{4.15}
S_{2} \le C \biggl(\sum\limits_{k=0}^{\infty}\biggl\{\sum\limits_{j=k}^{\infty}\biggl(\sum\limits_{m \in \mathbb{Z}^{n}}t^{p}_{k,m}2^{knp}\|v^{l}_{j}  |L_{1}(cQ^{n}_{k,m})\|^{p}\biggr\}^{\frac{\mu}{p}}\biggr)^{\frac{q}{\mu}}\biggr)^{\frac{1}{q}}.
\end{equation}

We next use inequality \eqref{eq4.4}, apply H\"older's inequality with exponents $p_{r}$, $p'_{r}$ to the sum over~$m$, and take into account Remark~\ref{R2.1}. As a~result,
\begin{gather}
\sum\limits_{m \in \mathbb{Z}^{n}}t^{p}_{k,m}2^{knp}\|v^{l}_{j}  |L_{1}(cQ^{n}_{k,m})\|^{p} \le C \sum\limits_{m \in \mathbb{Z}^{n}}t^{p}_{k,m}
\biggl(2^{kn}\sum\limits_{\substack{\widetilde{m} \in \mathbb{Z}^{n}\\
Q^{n}_{j,\widetilde{m}} \subset cQ^{n}_{k,m}}}\|v^{l}_{j}  |L_{1}(Q^{n}_{j,\widetilde{m}})\|\biggr)^{p} \le
\notag  \\
\le C \sum\limits_{m \in \mathbb{Z}^{n}}t^{p}_{k,m}2^{\frac{knp}{r}}\biggl(\sum\limits_{\substack{\widetilde{m} \in \mathbb{Z}^{n}\\
Q^{n}_{j,\widetilde{m}} \subset cQ^{n}_{k,m}}}2^{jn}\int\limits_{Q^{n}_{j,\widetilde{m}}}\frac{t^{r}_{j}(x)}{t^{r}_{j}(x)}\|v^{l}_{j}  |L_{r}(Q^{n}_{j,\widetilde{m}})\|^{r}\biggr)^{\frac{p}{r}} \le
\notag \\
\le C \sum\limits_{\widetilde{m} \in \mathbb{Z}^{n}}\biggl(\frac{\alpha^{1}_{k}}{\alpha^{1}_{j}}\biggr)^{p}t^{p}_{j,\widetilde{m}}2^{\frac{jnp}{r}}\|v^{l}_{j}  |L_{r}(cQ^{n}_{j,\widetilde{m}})\|^{p}.
\label{4.16}
\end{gather}

Substituting \eqref{4.16} into \eqref{4.15} and using Theorem~\ref{Th2.1}, this gives
\begin{equation}
\label{4.17}
\begin{split}
S_{2} \le C \|f|\widetilde{B}^{l}_{p,q,r}(\mathbb{R}^{n},\{t_{k}\})\|.
\end{split}
\end{equation}

To estimate $S_{1}$, we first write down the following chain of inequalities for $k \geq j$
(this follows from \eqref{2.5} and since
the overlapping multiplicity of the supports of splines $N^{l}_{j,m}$ is finite)
\begin{gather}
\allowdisplaybreaks
\int t^{p}_{k}(x)|\Omega_{k} \ast v^{l}_{j}|^{p}\,dx \le \notag \\
\le C \int t^{p}_{k}(x)\biggl[2^{(k-1)n}\int\limits_{\frac{I^{n}}{2^{(k-1)}}}|\Delta^{l}(h)v^{l}_{j}(x)|\,dh+2^{kn}\int\limits_{\frac{I^{n}}{2^{k}}}|\Delta^{l}(h)v^{l}_{j}(x)|\,dh\biggr]^{p} \le
\notag \\
\le C 2^{(j-k)lp} \int t^{p}_{k}(x)\biggl[\sum\limits_{\substack{ m \in \mathbb{Z}^{n}\\
x \in \operatorname{supp}N^{l}_{j,m}}}|\beta_{j,m}|\biggr]\,dx \le
\notag \\
\le C 2^{(j-k)lp}\sum\limits_{m \in \mathbb{Z}^{n}} 2^{\frac{njp}{r}}\biggl[\sum\limits_{\substack{\widetilde{m} \in \mathbb{Z}^{n}\\
Q^{n}_{k,\widetilde{m}} \subset cQ^{n}_{j,m}}} t^{p}_{k,\widetilde{m}}\biggr]\|v^{l}_{j} |L_{r}(cQ^{n}_{j,m})\|^{p}.
\label{4.18}
\end{gather}

Using \eqref{4.18} we have in view of (2.4)
\begin{gather}
\biggl(\int t^{p}_{k}(x) |\Omega_{k} \ast V^{l}_{1,k}(x)|^{p} \,dx\biggr)^{\frac{1}{p}}
\le \biggl(\int t^{p}_{k}(x)\biggl( \sum \limits_{j=0}^{k}  |\Omega_{k} \ast v^{l}_{j}(x)|^{\mu}\biggr)^{\frac{p}{\mu}} \,dx\biggr)^{\frac{1}{p}} \le
\notag \\
\le \biggl(\biggl(\frac{\alpha^{2}_{k}}{2^{kl}}\biggr)^{\mu}\sum \limits_{j=0}^{k}\biggl(\frac{2^{jl}}{\alpha^{2}_{j}}\biggr)^{\mu}
\biggl(\int t^{p}_{k}(x)|\Omega_{k} \ast v^{l}_{j}(x)|^{p}\biggr)^{\frac{\mu}{p}}\biggr)^{\frac{1}{\mu}} \le
\notag \\
\le C \biggl(\biggl(\frac{\alpha^{2}_{k}}{2^{kl}}\biggr)^{\mu} \sum \limits_{j=0}^{k} \biggl(\frac{2^{jl}}{\alpha^{2}_{j}}\biggr)^{\mu}
\biggl( \sum\limits_{\widetilde{m} \in \mathbb{Z}^{n}}t^{p}_{j,\widetilde{m}}2^{\frac{njp}{r}}\|v^{l}_{j}  |L_{r}(cQ^{n}_{j,m})\|^{p}\biggr)^{\frac{\mu}{p}}\biggr)^{\frac{1}{\mu}}.
\label{4.19}
\end{gather}

Applying estimate \eqref{4.19} in combination with Theorem~\ref{Th2.1} we obtain
\begin{equation}
\label{4.20}
S_{1} \le C \|f|\widetilde{B}^{l}_{p,q,r}(\mathbb{R}^{n},\{t_{k}\})\|.
\end{equation}

Now the required estimate follows from estimates \eqref{4.14}, \eqref{4.17}, \eqref{4.20}.

This completes the proof of the theorem.

We are now ready to formulate the main result of the paper.

\begin{corollary}\label{Ca}
Let $p,q,r \in (0,\infty]$, $\theta \in (0,\min\{1,r,p\}]$.  Next, let a function $\varphi_{0} \in D(\mathbb{R}^{n})$,
$\varphi=\varphi_{0}-2^{-n}\varphi_{0}(\frac{\cdot}{2})$, and let $p=\{t_{k}\}$~be an admissible
weight sequence,  $\{t_{k}\} \in X^{\alpha_{3}}_{\alpha,\sigma,p}$ with $\sigma_{1}=\theta p'_{\theta}$, $\sigma_{2}=p$, $\{\alpha^{1}_{k}\}=\{2^{k\alpha_{1}}\}$,
$\{\alpha^{2}_{k}\}=\{2^{k\alpha_{2}}\}$, $\min\{l,L_{\varphi}+1\} > \alpha_{2}$, $\alpha_{1} > n(\frac{1}{\theta}-\frac{1}{\max\{r,1\}})$. Then
$B^{\varphi_{0}}_{p,q}(\mathbb{R}^{n}, \{t_{k}\})=\widetilde{B}^{l}_{p,q,r}(\mathbb{R}^{n},\{t_{k}\})$, the corresponding norms being equivalent.
\end{corollary}

\textbf{Proof.} From Lemma \ref{Lm2.1} it follows that $\{\overline{t}_{k}\} \in X^{\alpha_{3}}_{\alpha,\sigma,p}$. Since $L_{\varphi}+1 > \alpha_{2}$ we have from
Theorem~\ref{Th3.3} that $B^{\varphi_{0}}_{p,q}(\mathbb{R}^{n},\{t_{k}\})=B^{\varphi_{0}}_{p,q}(\mathbb{R}^{n},\{\overline{t}_{k}\})$,
the corresponding norms being equivalent. Hence, in view of the condition $\alpha_{1} > n(\frac{1}{\theta}-\frac{1}{\max\{r,1\}})$
an application of Theorem~\ref{Th4.1} shows that  $B^{\varphi_{0}}_{p,q}(\mathbb{R}^{n},\{t_{k}\}) \subset L^{\text{\rm loc}}_{1}(\mathbb{R}^{n})$.
Hence, combining Theorems \ref{Th4.2}, \ref{Th4.3} (the assumptions of these theorems are clearly satisfied) this gives $\widetilde{B}^{l}_{p,q,r}(\mathbb{R}^{n},\{t_{k}\})=\widetilde{B}^{l}_{p,q,\theta}(\mathbb{R}^{n},\{t_{k}\})=B^{\varphi_{0}}_{p,q}(\mathbb{R}^{n},\{t_{k}\})$,
the corresponding norms are equivalent.

\begin{Remark} {\rm
Let us consider some particular cases of Corollary~\ref{Ca}.

Applying Corollary \ref{Ca} with $r=1$, $\theta = \min\{r,p\}$, taking into account Remarks \ref{R2.3}, \ref{R2.4}, we obtain in view of Theorem 2.5 of~\cite{Tyu}
that, for $p,q \in (0,\infty]$, $\{s_{k}\} \in Y^{\alpha_{3}}_{\alpha_{1},\alpha_{2}}$, $\alpha_{1} > n\max\{\frac{1}{p}-1,0\}$, $l > \alpha_{2}$,
the space $B^{\{s_{k}\}}_{p,q}(\mathbb{R}^{n})$ consists of the functions $f \in L^{\text{\rm loc}}_{\max\{1,p\}}(\mathbb{R}^{n})$ (moreover,
$B^{\{s_{k}\}}_{p,q}(\mathbb{R}^{n})$ is continuously embedded  into $L^{\text{\rm loc}}_{\max\{1,p\}}(\mathbb{R}^{n})$). Besides,
\begin{equation}
\|f|B^{\{s_{k}\}}_{p,q}(\mathbb{R}^{n})\| \sim \left(\sum\limits_{k=1}^{\infty}\|s_{k}\overline{\Delta}^{l}_{1}(2^{-k})f|L_{p}(\mathbb{R}^{n})\|^{q}\right)^{\frac{1}{q}}+\|s_{0}\|f|L_{1}(\cdot+I^{n})\||L_{p}(\mathbb{R}^{n})\|
\end{equation}
(the modification in the case $q=\infty$ are standard), where, for a  function $f \in L^{loc}_{r}(\mathbb{R}^{n})$ ($r \in (0,\infty)$), $l\in \mathbb N$  we set

$$
\overline{\Delta}^{l}_{r}(2^{-k})f(x):=\Bigl(2^{kn}\int\limits_{\frac{I^{n}}{2^{k}}}|\Delta^{l}(h)f(x)|^{r}\,dh \Bigr)^{\frac{1}{r}}, \qquad x \in \mathbb{R}^{n},k \in \mathbb{N}.
$$

This result is slightly different from Theorem 18 of~\cite{Kempka}. Indeed, as distinct from~\cite{Kempka},
we do not impose any restrictions on the parameter $\alpha_{3}$, but at the same time we do not guarantee that the embedding $B^{\{s_{k}\}}_{p,q}(\mathbb{R}^{n}) \subset L_{\max\{1,p\}}(\mathbb{R}^{n})$ is continuous.

Assume now that $p,q,r \in (0,\infty]$, $p,r \neq \infty$, $\theta \in (0,p]$, $\gamma^{p} \in A_{\frac{p}{\theta}}(\mathbb{R}^{n})$.

Taking into account Theorem 2.5 of \cite{Tyu} and Remarks \ref{R2.3}, \ref{R2.4} we see that, for $s > n(\frac{1}{\theta}-\frac{1}{\max\{1,r\}})$, the
weighted Besove space  $B^{s}_{p,q}(\mathbb{R}^{n},\gamma)$ is continuously embedded into
$L^{\text{\rm loc}}_{\max\{1,p\}}(\mathbb{R}^{n})$). Moreover,
\begin{equation}
\|f|B^{s}_{p,q}(\mathbb{R}^{n},\gamma)\| \sim \left(\sum\limits_{k=1}^{\infty}\|\gamma \overline{\Delta}^{l}_{r}(2^{-k})f|L_{p}(\mathbb{R}^{n})\|^{q}\right)^{\frac{1}{q}}+\|\gamma\|f|L_{r}(\cdot+I^{n})\||L_{p}(\mathbb{R}^{n})\|
\end{equation}
(the modification for $q=\infty$ is standard).

Thus result is a~particular case of Theorem 3.14 of~\cite{HN} and was obtained using different technique.
}
\end{Remark}

\section{Applications}
As a possible application of the above results we give an equivalent description of the trace space of the weighted Sobolev space.

Throughout this section we fix natural numbers $d$, $n > d$ and a~parameter $p \in (1,\infty)$.
A~point of the space $\mathbb{R}^{n+d}=\mathbb{R}^{n} \times \mathbb{R}^{d}$ will be written as the pair $x=(x',x'')$.

Let $l \in \mathbb{N}$ and let $\gamma$ be a~weight. By $W^{l}_{p}(\mathbb{R}^{n},\gamma)$ we shall denote the weighted Sobolev space with the norm
$$
\|f|W^{l}_{p}(\mathbb{R}^{n},\gamma)\|:=\sum\limits_{\alpha \le l}\|\gamma D^{\alpha}f|L_{p}(\mathbb{R}^{n})\|,
$$
where $D^{\alpha}f$ are the (Sobolev) generalized derivatives of a~function~$f$ of order~$\alpha$.

In what follows, we shall assume that the weight lies in the weighted Muckenhoupt class $A_{p}(\mathbb{R}^{n})$ (the definition of the class $A_{p}(\mathbb{R})$
may be found, for example, in Ch.~5 of~\cite{St}).

Given $k \in \mathbb{N}_{0}$, we set
$$
\gamma^{p}_{k}(x'):=\sum\limits_{m \in \mathbb{Z}^{n}}\chi_{\widetilde{Q}^{n}_{k,m}}(x')\iint\limits_{Q^{n}_{k,m}\times (\frac{B^{d}}{2^{k}} \setminus \frac{B^{d}}{2^{k+1}})}\gamma(x',x'')\,dx'dx'',  x' \in \mathbb{R}^{n}.
$$

It is easily seen that $\{\gamma_k\}\in Y^{\alpha_3}_{\alpha_1,\alpha_2}$ for some $\alpha_1, \alpha_2\in \mathbb R$, $\alpha_3\ge 0$.

We do not define the trace on the plane of a~weighted Sobolev space, but refer the reader to~\cite{Tyu2}.

The following result is obtained by a~combination  of Theorem 3.1 of \cite{Tyu} with Corollary~\ref{Ca} and Remark \ref{R2.4}.

\begin{theorem}
\label{Sled}
Let  $p \in (1,\infty)$, $r \in [1,p) $,  $\gamma^{p} \in
A_{\frac{p}{r}}(\mathbb{R}^{n})$, $f \in
W_{p}^{l}(\mathbb{R}^{n},\gamma)$, $l > \frac{n-d}{r}$. Then there exists the trace $\varphi \in
B^{\{\gamma_{k}\}}_{p,p}(\mathbb{R}^{d})$ of the function~$f$, and moreover,
\begin{equation}
\label{trineq}
\left\|\varphi|B^{\{\gamma_{k}\}}_{p,p}(\mathbb{R}^{d})\right\|
\le C_{1}\|f|W_{p}^{l}(\mathbb{R}^{n},\gamma)\|.
\end{equation}
The constant $C_{1}$ in~\eqref{trineq} is independent of the function~$f$.

Conversely, if a function $\varphi \in
B^{\{\gamma_{k}\}}_{p,p}(\mathbb{R}^{d})$, then there exists a~function $f \in
W_{p}^{l}(\mathbb{R}^{n},\gamma)$ such that $\varphi$ is the trace of~$f$ on $\mathbb{R}^{d}$, and moreover,
\begin{equation}
\label{extineq}
\|f\,|\,W_{p}^{l}(\mathbb{R}^{n},\gamma)\|
\le C_{2} \bigl\|\varphi \,\big|\,B^{\{\gamma_{k}\}}_{p,p}(\mathbb{R}^{d})\bigr\|,
\end{equation}
the constant $C_{2}$ in \eqref{extineq} being independent of the function~$\varphi$.
\end{theorem}

\begin{Remark}{\rm
It is worth noting that this result is new and cannot be obtained by the methods available before.
Indeed, the previously available machinery, which was developed for studying Besov spaces of variable smoothness, is much like the classical methods
for dealing with Besov spaces of constant smoothness. The methods known so far depend on pointwise estimates of the variable smoothness, and hence
are incapable of yielding Theorem~5.1 (see \cite{Tyu} for details).}
\end{Remark}



\section*{References}

\begin{thebibliography}{10}

\bibitem{Almeida}
A.~Almeida and P.~H{\"a}st{\"o},
Besov spaces with variable smoothness and integrability,
\textit{J. Funct. Anal}. \textbf{258}:5, 1628--1655 (2010).

\bibitem{Fix}
 C.~de Boor and G.F.~Fix, Spline approximation by quasi-interpolants,\textit{J. Approx. Theory}, 8, 19--45 (1973).

\bibitem{Bennet}
 G.~ Bennet, Some elementary inequalities, \textit{Quart. J. Math.Oxford} Ser. (2), 38, 401--425 (1987).

\bibitem{Be1}
 O.\,V.~Besov, Equivalent normings of spaces of functions of variable smoothness, \textit{Proceedings of the Steklov Institute of Mathematics}, \textbf{243}, 80--88 (2003).

\bibitem{Be2}
 O.\,V.~Besov, Interpolation, embedding, and extension of spaces of functions of variable smoothness, \textit{Proc. Steklov Inst. Math.} \textbf{248}, 47--58 (2005);
translation from \textit{Tr. Mat. Inst. Steklova} \textbf{248}, 52--63 (2005).

\bibitem{DeVore}
 R.A.~ DeVore, V.A.~ Popov, Interpolation of Besov spaces, \textit{Trans. Amer. Math. Soc.} \textbf{305} (1) (1988) 397--414.

\bibitem{HN}
 L.I.~Hedberg and Y.~Netrusov, An axiomatic approach to function
spaces, spectral synthesis, and Luzin approximation, \textit{Mem. Amer. Math. Soc}. \textbf{188} (882) (2007),

\bibitem{Kempka2}
 H. Kempka, Atomic, molecular and wavelet decomposition of generalized 2-microlocal Besov spaces, \textit{J. Funct. Spaces Appl.}, \textbf{8} (2010), no.~2, 129--165.

\bibitem{Kempka}
 H.~Kempka and J.~Vybiral, Spaces of variable smoothness and
integrability: characterizations by local means and ball means of differences. \textit{J.~Fourier Anal. Appl.} \textbf{18} (2012), no. 4, 852--891.

\bibitem{Kempka}
 H.~Kempka, 2-microlocal Besov and Triebel-Lizorkin spaces of variable integrability, \textit{Revista Matem\'atica Complutense},
 \textbf{22}, 1, 227--251 (2009).
 

\bibitem{Dachun}
Y.~Liang, D. Yang, W. Yuan, Y. Sawano, and T. Ullrich, A new framework for generalized Besov-type and Triebel--Lizorkin-type spaces, \textit{Diss. Math.}
(\textit{Rozprawy Mat.}) \textbf{489}, 114(2013).

\bibitem{Rudin}
W.~Rudin, Functional Analysis, \textit{McGraw-Hill}, New York, 1991.


\bibitem{Ry}
 V.~S.~Rychkov, Littlewood--Paley theory and function spaces with $A_{p}^{\text{\rm loc}}$-weights, \textit{Math. Nachr.}, \textbf{224}, 145--180 (2001).

\bibitem{Ry2}
V.~S.~Rychkov, On restrictions and extensions of the Besov and Triebel-Lizorkin spaces with respect to Lipschitz domains, J. London Math. Soc. (2), \textbf{60}, 1, 237--257 (1999).


\bibitem{Tri}
H.~Triebel, Theory of Function Spaces, \textit{Birkh{\"a}user},
Basel,1983.


\bibitem{Tyu}
A.I.~Tyulenev, Some new function spaces of variable smoothness,
 \textit{Sbornik: Mathematics}, \textbf{6}  (2015) to appear. Preliminary version http://arxiv.org/abs/1410.8360

\bibitem{Tyu2}
 A.~I.~Tyulenev, Description of traces of functions in the Sobolev space with a Muckenhoupt weight, \textit{Proc. Stekl. Inst. Math.} \textbf{284} (2014), 280--295.

\bibitem{St}
E.~M.~Stein, Harmonic Analysis: Real-Variable methods, Orthogonality, and Oscillatory Integrals, \textit{Princeton Univ. Press}, Princeton, NJ, 1993.

\bibitem{Ullrich}
T.~ Ullrich, H.~ Rauhut: Generalized coorbit space theory and inhomogeneous function spaces of Besov--Lizorkin--Triebel type, \textit{J. Funct. Anal.} \textbf{260} (11), 3299--3362 (2011).


\end{thebibliography}

\end{document}